\title{\normalsize \large
 {\textbf     {
THE SIGN CLUSTERS OF THE MASSLESS GAUSSIAN  \\ FREE FIELD PERCOLATE ON $\mathbb{Z}^d$, $d\geqslant3$ (AND MORE) 
 }}}

\date{}
\author{}
\documentclass[12pt, leqno]{article}
\usepackage[english]{babel}
\usepackage[
left=2.7cm,right=2.7cm,top=2.7cm,bottom=3.5cm]{geometry}
\usepackage{amsfonts, amsmath, amssymb, amsthm, mathrsfs}
\usepackage{graphicx}
\usepackage{color}
\usepackage[sans]{dsfont}
\usepackage{url}

\usepackage{color}

\usepackage{comment}
\usepackage[hang,flushmargin]{footmisc}

\numberwithin{equation}{section}

\newtheorem{The}{Theorem}[section]
\newtheorem{Lemme}[The]{Lemma}
\newtheorem{Prop}[The]{Proposition}
\newtheorem{Cor}[The]{Corollary}
\theoremstyle{definition}
\newtheorem{Def}[The]{Definition}

\newtheorem{Rk}[The]{Remark}

\theoremstyle{remark}

\theoremstyle{definition}

\usepackage{csquotes}
\usepackage{IEEEtrantools}
\usepackage{mathtools}
\usepackage{enumerate}

\makeatletter
\renewcommand{\@maketitle}{
\newpage
 \null
 \vskip 2em%
 \begin{center}%
  {\LARGE \@title \par}%
 \end{center}%
 \par} \makeatother

\makeatletter
\newcommand{\E}{\mathbb{E}}

\newcommand{\Q}{\mathbb{Q}}
\newcommand{\R}{\mathbb{R}}
\newcommand{\Z}{\mathbb{Z}}
\newcommand{\N}{\mathbb{N}}
\renewcommand{\P}{\mathbb{P}}
\newcommand{\eps}{\varepsilon}

\newcommand{\tp}{\widetilde{\varphi}}
\newcommand{\I}{{\cal I}}
\newcommand{\V}{{\cal V}}
\renewcommand{\phi}{\varphi}
\renewcommand{\tilde}{\widetilde}
\renewcommand{\hat}{\widehat}
\renewcommand{\epsilon}{\varepsilon}
\renewcommand\theequation{\thesection.\arabic{equation}}

\@addtoreset{equation}{section}
 
\makeatother
\usepackage{color}

\begin{document}

\maketitle

\begin{center}
\vspace{0.7cm}
Alexander Drewitz$^1$, Alexis Pr\'evost$^1$ and Pierre-Fran\c cois Rodriguez$^2$

\end{center}
\vspace{0.1cm}
\begin{abstract}
\centering
\begin{minipage}{0.8\textwidth}
\vspace{0.3cm}
We investigate the percolation phase transition for level sets of the Gaussian free field on $\mathbb{Z}^d$, with $d\geqslant 3$, and prove that the corresponding critical parameter $h_*(d)$ is strictly positive for all $d\geqslant3$, thus settling an open question 
from \cite{MR3053773}. 
In particular, this implies that the sign clusters of the Gaussian free field percolate on $\mathbb{Z}^d$, for all $d\geqslant 3$. Among other things, our construction of an infinite cluster above small, but positive level $h$ involves random interlacements at level $u>0$, a random subset of $\mathbb{Z}^d$ with desirable percolative properties, introduced in \cite{MR2680403} in a rather different context, a certain Dynkin-type isomorphism theorem relating random interlacements to the Gaussian free field \cite{MR2892408}, and a recent coupling of these two objects \cite{MR3502602}, lifted to a continuous metric graph structure over~$\mathbb{Z}^d$.
\end{minipage}
\end{abstract}

\thispagestyle{empty}

\vspace{2.5cm}

\begin{flushleft}

\noindent\rule{5cm}{0.4pt} \hfill May 2018 \\
\bigskip
$^1$Universit\"at zu K\"oln\\
Mathematisches Institut \\
Weyertal 86--90 \\
50931 K\"oln, Germany. \\
\url{drewitz@math.uni-koeln.de}\\
\url{aprevost@math.uni-koeln.de} \\

\bigskip
\medskip
$^2$Department of Mathematics \\
University of California, Los Angeles \\
520, Portola Plaza, MS 6172\\
Los Angeles, CA 90095 \\
\texttt{rodriguez@math.ucla.edu}
\end{flushleft}

\newpage

\section{Introduction}

The present work studies the percolation phase transition of Gaussian free field level sets on $\Z^d,\ d\geqslant3$, which provides a canonical example for a percolation model with strong, algebraically decaying correlations. It was first proved in \cite{MR914444} that the corresponding critical level $h_*(d)$, see \eqref{eq:h_*} below for its definition, satisfies $h_*(d) \geqslant 0$ for every dimension $d \geqslant3$ and that $h_*(3) < \infty$. It was later shown in \cite{MR3053773} that $h_*(d)$ is finite in every dimension $d \geqslant3$, and strictly positive when $d$ is large, with leading asymptotics as $d \to \infty$ derived in \cite{MR3339867}. We prove here that this parameter is actually strictly positive in all dimensions $ d \geqslant3$. This answers a question from \cite{MR914444}, see also Remark 3.6 in \cite{MR3053773}, and fits with numerical evidence from \cite{Marinov}, see Section 4.1.2 and Figure 4.1 therein. A corresponding classical result for Bernoulli site percolation, $p_c^{\text{site}}(\mathbb{Z}^d)<\frac{1}{2}$ for $d \geqslant 3$, has been known to hold for several decades already \cite{MR781418}.

Our construction of infinite clusters (by which, adopting the usual terminology, we mean unbounded connected components) of excursion sets for the Gaussian free field crucially relies on another object, random interlacements. The model of random interlacements has originally been introduced in \cite{MR2680403} to study certain geometric properties of random walk trajectories on large, asymptotically transient, finite graphs. The Dynkin-type isomorphism theorem relating interlacements and the Gaussian free field, see \cite{MR2892408}, has repeatedly proved a useful tool in their study, see \cite{MR2892408}, \cite{MR2932978}, \cite{MR3163210}, \cite{MR3502602}, \cite{MR3492939} and \cite{AbacherliandSznitman}. In a broader scheme, the usefulness of similar random walk representations as a tool in field theory and statistical mechanics has been recognized for a long time, see \cite{Sym}, \cite{brydges1982random} and \cite{dynkin1983markov}.

The cable system method introduced in \cite{MR3502602} provides a continuous version of this isomorphism theorem, from which some links between the level sets of the Gaussian free field and the vacant sets $\V^u$, $u > 0$, of random interlacements can be derived. This method was used in \cite{MR3492939} and \cite{AbacherliandSznitman} to find a suitable coupling between those two sets, and was applied in the case of transient trees. It was also proved in these papers that, under certain conditions on the geometry of the tree $\mathbb{T}$, the critical parameter $h_*(\mathbb{T})$ for level set percolation of the Gaussian free field on $\mathbb{T}$ is strictly positive. As will become apparent below, the isomorphism theorem on the cable system can be paired with renormalization techniques from random interlacements, and in particular from \cite{MR3024098}, which imply a certain robustness property of $\mathcal{I}^u = \Z^d \setminus \V^u $ with respect to small noise, to yield similar findings on $\Z^d$, for all $d\geqslant3$.

\medskip
Let us now describe the results in more details. For $d\geqslant3,$ we consider $\Z^d$ as a graph, with undirected edge set $E$, and take uniform weights equal to 1 on all edges in $E$, so  that the sum of all weights around a vertex $x\in{\Z^d}$ is $2d$. For $x,y\in{\Z^d},$ we write $x\sim y$ if and only if $\{x,y\}\in{E}$. Noting that $\Z^d$, $d\geqslant3$, is transient for discrete time simple random walk, we define the symmetric Green function by
\begin{equation}
\label{green}
    g(x,y)=\frac{1}{2d} E_x\Big[\int_{0}^\infty1_{\{X_t=y\}} \, \mathrm{d}t\Big], \ x,y\in{\Z^d},
\end{equation}
where $(X_t)_{t\geq 0}$ denotes the canonical continuous time random walk on $\Z^d,$ with constant jump rate $1,$ starting at $x$ under $P_x.$ We also set $g(x)=g(0,x)$, for $x \in \mathbb{Z}^d$. We define $\P^G$, a probability measure on $\R^{\Z^d}$ endowed with its canonical $\sigma$-algebra generated by the coordinate maps $\Phi_x$, $x\in{\Z^d}$, such that, under~$\P^G$,
\begin{equation}
\begin{split}
\label{Phi}
    &(\Phi_x)_{x\in{\Z^d}}\text{ is a centered Gaussian field with}\\
    &\text{covariance function } \E^G[\Phi_x\Phi_y]=g(x,y)\text{ for all }x,y\in{\Z^d}.
    \end{split}
\end{equation}
(Any random field $\phi=(\phi_x)_{x\in{\Z^d}}$ with law $\P^G$ on $\R^{\Z^d}$ will henceforth be called a Gaussian free field on ${\Z}^d$). We are interested in level sets of $\Phi$, and for every $h\in{\R},$ denote by $\{x\stackrel{\mathrm{\geqslant}h}{\longleftrightarrow}\infty\}$ the event that $x\in{\Z^d}$ lies in an infinite connected component of 
\begin{equation}
\label{exc}
E^{\geqslant h} = \{y\in{\Z^d};\, \Phi_y\geqslant h\},
\end{equation} 
and by $\eta(h)$ its probability, which does not depend on the choice of $x.$ 
 The function $\eta(\cdot)$ is decreasing, and it is natural to ask whether it is strictly positive or not. This leads to the definition of the critical point \begin{equation}
\label{eq:h_*}
    h_*(d) \stackrel{\text{def.}}{=} \inf\left\{h\in{\R};\, \eta(h)=0\right\}.
\end{equation}
By ergodicity, this definition corresponds to the phase transition for the existence of an infinite connected component in $E^{\geqslant h}$, see Lemma 1.5 in \cite{MR3053773}. It is not a priori clear whether $|h_*|<\infty$ or not, and a summary of the \textit{status quo} was given in the first paragraph. In summary, it is known that $h_*(d) \in [0,\infty)$ for all $d\geqslant 3$, and that $h_*(d)\sim (2g(0)\log d)^{1/2}$ as $d \to \infty$. Our main result is the following lower bound in all dimensions.
\begin{The}
\label{mainresult}
\begin{equation}
\label{eq:mainresult}
    h_*(d)> 0, \quad \text{for all $d \geqslant 3$.}
\end{equation}
Moreover, there exists $h_1>0$ such that for all $h\leqslant h_1,$ there exists $L_0>0$ such that $E^{\geqslant h}$ contains an infinite cluster in the thick slab $\Z^2\times[0,2L_0)^{d-2}.$

\end{The}

In fact, one can replace $E^{\geqslant h} $ by $\{y\in{\Z^d};\,   K(h) \geqslant \Phi_y\geqslant h \} (\subset E^{\geqslant h} ) $ for sufficiently large $K(h)$ in the previous statement, see Remark \ref{R:final}, 2) below. Note that the infinite cluster of $E^{\geqslant h}$ cannot be contained in $\Z^2\times\{0\}^{d-2}$ for $0\leqslant h<h_*(d)$, as explained in Remark \nolinebreak3.6.1 of \cite{MR3053773}. As an immediate corollary of Theorem \ref{mainresult}, we note that there exists an open interval $I \subset \mathbb{R}$ containing the origin and such that, for all $h \in {I}$, the level set $E^{\geqslant h}$ and its complement $E^{< h} = \mathbb{Z}^d \setminus E^{\geqslant h}$ \textit{both} percolate (with probability one). This follows readily from  \eqref{eq:mainresult} and the fact that $E^{< h} \stackrel{\text{law}}{=} E^{\geqslant - h}$ for all $h \in \mathbb{R}$, by symmetry of $\Phi$. In particular, choosing $h=0$, this implies that 
\begin{equation}
\label{eq:1.10}
\text{$\Phi$ almost surely contains two infinite sign clusters (one for each sign).}
\end{equation}
Put differently, Theorem \ref{mainresult} asserts that the critical density $p_c^G (d) = \mathbb{P}^G[\Phi_0 \geqslant h_*(d)]$ satisfies $p_c^G (d) < \frac12$, for all $d \geq 3$, thus mirroring the result $p_c^{\text{site}}(\mathbb{Z}^d)<\frac{1}{2}$, see  \cite{MR781418}, for independent Bernoulli site percolation on $\mathbb{Z}^d$, $d\geq 3$. However, the elegant geometric arguments developed therein to ``interpolate'' between two- and three-dimensional structures do not seem to transfer to the current situation: the correlations present a serious impediment. Moreover, there is no obvious monotonicity of $p_c^G (d)$ (or $h_*(d)$) as a function of $d$. One may also conjecture that $p_c^G(d) < p_c^{\text{site}}(\Z^d)$, the critical density for independent site percolation on the lattice, based on the reasonable intuition that positive correlations ``help'' in forming clusters of $E^{\geqslant h}$. We do not currently know a proof of this (nor of the more modest conjecture $p_c^G(d) \leq p_c^{\text{site}}(\Z^d)$).


\medskip

A key tool in the proof of Theorem \ref{mainresult} is a certain isomorphism, see Theorem \ref{iso} below, which gives a link between random interlacements and the Gaussian free field. We now explain its benefits in some detail, and refer to Section \ref{cablesystem} for precise definitions. Suppose that $\omega$ denotes the interlacement point process defined in \cite{MR2680403}, with law $\P^I$, and let $\omega^u$ be the process consisting of the trajectories in the support of $\omega$ with label at most $u$. Somewhat informally, $\omega^u$ is a Poisson cloud of bi-infinite nearest neighbor trajectories modulo time-shift whose forward and backward parts escape all finite sets in finite time. One naturally associates to $\omega^u$, see for instance (1.8) in \cite{MR2892408}, a field of occupation times $(\ell_{x,u})_{x\in{\Z^d}}$, where $\ell_{x,u}= \ell_{x,u}(\omega^u)$ collects the total amount of time spent at $x$ by any of the trajectories in the support of $\omega^u$. The interlacement set at level $u$ is then defined as
\begin{equation}
\label{Iu}
    \I^u=\{x\in{\Z^d};\, \ell_{x,u}>0\}.
\end{equation}
It corresponds to the set of vertices visited by at least one trajectory in the support of $\omega^u$. For any $u > 0$, the set $\I^u$ is almost surely unbounded and connected \cite{MR2680403}. The following isomorphism was proved in Theorem 0.1 of \cite{MR2892408}, and has the same spirit as the generalized second Ray-Knight theorem, see for example \cite{MR1813843}, \cite{MR2250510} or \cite{MR2932978}:
\begin{equation}
\label{discreetiso}
    \left(\ell_{x,u}+\frac{1}{2}\Phi_x^2\right)_{x\in{\Z^d}} \text{ under $\P^I\otimes\P^G$ has the same law as } \left(\frac{1}{2}(\Phi_x+\sqrt{2u})^2\right)_{x\in{\Z^d}}\text{ under }\P^G.
\end{equation}

If one attaches to each edge $e$ of $\Z^d$ a line segment $I_e$ of length $\frac{1}{2},$ the resulting ``graph'' $\tilde{\Z}^d$ is continuous and called the \textit{cable system}, see Section \ref{cablesystem}. On this cable system, one then defines probabilities $\tilde{\P}^G$ and $\tilde{\P}^I$ under which the fields $(\Phi_x)_{x\in\mathbb{Z}^d}$ and $(\ell_{x,u})_{x\in \Z^d}$ admit \textit{continuous} extensions $\widetilde{\Phi}= (\tilde{\Phi}_x)_{x\in{\tilde{\Z}^d}}$ and $\tilde{\ell}=(\tilde{\ell}_{x,u})_{x\in{\tilde{\Z}^d}},$ and the set $\tilde{\I}^u=\{x\in{\tilde{\Z}^d};\,\tilde{\ell}_{x,u}>0\}$ is connected. It was proved in \cite{MR3502602} that for each $u>0,$ a continuous version of the isomorphism \eqref{discreetiso} also holds on $\tilde{\Z}^d$, see also \eqref{isoeq} below, and in particular (somewhat inaccurately, but see \eqref{isoeq}, \eqref{htilde=0} below for precise statements) the sign of $\tilde{\Phi}_x+\sqrt{2u}$ is constant as long as $\tilde{\ell}_{x,u}>0,$ and thus by the continuity of $\tilde{\Phi}$ and the connectivity of $\tilde{\I}^u$, either $\tilde{\Phi}_x>-\sqrt{2u}$ for all $x\in{\tilde{\I}^u}$ or $\tilde{\Phi}_x<-\sqrt{2u}$  for all $x\in{\tilde{\I}^u}.$ But $\tilde{\I}^u$ is unbounded, hence, taking $h=\sqrt{2u},$ by symmetry of the Gaussian free field and ergodicity, $\tilde{\P}^G$-a.s.\ the set
\begin{equation}
\label{unboun}
    \{x\in{\tilde{\Z}^d};\, \tilde{\Phi}_x\geqslant-h\}\text{ contains an unbounded cluster in the cable system $\tilde{\Z}^d.$}
\end{equation}
This result was already known to hold on $\Z^d$ without the isomorphism theorem  \cite{MR914444}, where it had been derived using a neat contour argument. It is interesting to note that, on the cable system, \eqref{unboun} is actually sharp, because $\tilde{\P}^G$-a.s.\ the set 
\begin{equation}
\label{eq:1.20}
    \{x\in{\tilde{\Z}^d}; \, \tilde{\Phi}_x\geqslant0\}\text{ does not contain unbounded clusters in the cable system $\tilde{\Z}^d,$}
\end{equation}
see Proposition 5.5 in \cite{MR3502602}, which sharply contrasts with \eqref{eq:1.10}. All in all, the infinite cluster in $E^{\geqslant0}$ (part of $\mathbb{Z}^d$), which exists by Theorem \ref{mainresult}, ``scatters'' into finite pieces upon adding the field on the edges, but the infinite cluster of $E^{\geqslant -h}$ does not, for ever so small $h>0$.


On our way towards proving Theorem \ref{mainresult}, we will first show that a truncated version of the level sets in \eqref{unboun} contains an unbounded cluster on $\tilde{\Z}^d.$ Indeed, it was proved in \cite{MR3024098} that the intersection of the random interlacement set $\I^u$ with a Bernoulli percolation having large success parameter still contains an infinite cluster in $\Z^d.$ By showing a similar stability result on the cable system, see Proposition \ref{truncated negative level sets}, and using the isomorphism theorem on the cable system, we will obtain, cf.\ Theorem \ref{percphitilde} below, that the truncated (continuous) level set 
\begin{equation}
\label{1:30}
\{x\in{\tilde{\Z}^d};\, -h\leqslant\tilde{\Phi}_{x}\leqslant K(h)\}
\end{equation}
contains an unbounded cluster on $\tilde{\Z}^d$ for all $h > 0$ and large enough, but finite $K(h)$ (with $hK(h)\to 0$ as $h \searrow 0$). Once this has been proved, see Theorem \nolinebreak\ref{percphitilde}  for the precise technical statement, we no longer need to use random interlacements to prove Theorem \ref{mainresult} (note however that the interlacements are crucial in generating a suitable percolating cluster to start with, i.e., one which is already reasonably ``close'' to being a sign cluster of the free field, see \eqref{1:30}). 

We now describe the second part of the proof. By construction, one can view $\tilde{\Phi}$, the Gaussian free field on the cable system~$\tilde{\Z}^d$, as a Gaussian free field on $\Z^d$ with Brownian bridges of length $\frac{1}{2}$ attached on the edges, see \eqref{bb} and thereafter. On an edge contained in the set of \eqref{1:30}, those Brownian bridges never go below $-h,$ which happens with low probability for small $h.$ We are going to use this low probability to go from $-h\leqslant\tilde{\Phi}\leqslant K(h)$
\textit{on the edges} to $h\leqslant\tilde{\Phi}\leqslant K(h)$ \textit{on the endpoints} of theses edges and for small enough $h$, see in particular Lemma \ref{comparingEandF}, which will then imply that the set $\{x\in{\Z^d}; \, \tilde{\Phi}_x\geqslant h\}$ has an infinite cluster on $\Z^d$, as asserted. 

\medskip

We now explain the organization of this article, and highlight its main contributions. In Section \ref{cablesystem}, we recall the definitions of the Gaussian free field and random interlacements on the cable system, and the link between the two via the aforementioned isomorphism theorem. In Section \ref{strong}, we collect a few preparatory tools by showing some strong connectivity properties, a large deviation inequality as well as a version of the decoupling inequalities for random interlacements on the cable system. Most of these are well-known in spirit, but existing results do not entirely fit our needs. 

The construction of the infinite cluster comes essentially in three steps, Proposition \ref{truncated negative level sets}, Theorem \ref{percphitilde}, and Section \ref{proofofthemainresult}, which are the main reference points of this paper. Proposition~\ref{truncated negative level sets} is a fairly generic result, which, roughly speaking, for \textit{any} coupling of a continuous interlacement and a Gaussian free field, see \eqref{eq:Qdef}, yields a percolating interlacement cluster, with good control on the free field part, and some room to play with along the edges. Its proof follows a standard static renormalization scheme from \cite{MR3024098}, \cite{MR3163210}, assembling the results of Section \ref{strong}. In Theorem \ref{percphitilde}, we ``translate'' Proposition \ref{truncated negative level sets}, for a certain choice of the coupling, to show that suitably truncated level sets of the Gaussian free field on the cable system contain an unbounded connected component. The reference level for the excursion sets of Theorem \ref{percphitilde} is $-h$, for (small) positive $h$. Section \ref{proofofthemainresult} contains the device to ``flip the sign'' and pass from $-h$ to $h$ on the vertices, as indicated above. Together with Theorem \ref{percphitilde}, this then yields a proof of Theorem \ref{mainresult}.

\medskip

In the rest of this article, we denote by $c$ and $C$ positive constants that may change from place to place. Numbered constants such as $C_0,$ $c_0,$ $C_1,$ $C'_1,\dots$ are fixed until the end of the article. All constants are allowed to implicitly depend on the dimension $d$ and a parameter $u_0>0$, which will first appear in Lemma \ref{connectivity} and throughout the remaining sections.

\section{Notation and useful facts about the cable system}
\label{cablesystem}
In this section, we give a definition of the Gaussian free field and random interlacements on the cable system that will be useful later. We also discuss some aspects of the Markov property for the Gaussian free field and its consequences, and recall the isomorphism theorem which links random interlacements and the Gaussian free field. 

For later convenience, we endow the graph $\Z^d$ with a distance function $d(\cdot,\cdot)$ which is half of the usual graph distance, i.e., half of the $\ell^1$-distance $|\cdot|_1$ on $\Z^d.$ Recall that we write $x \sim y$, for $x,y \in \mathbb{Z}^d$, if $|x-y|_1=1$. We define $V^0=\{2x,\ x\sim0\},$ so that, for all $x,y\in{\Z^d}$ with $x\sim y,$ we can write $y=x+\frac{1}{2}v_{(x,y)}$ for a unique $v_{(x,y)}\in{V^0}$. Note that $d(x,x+\frac{1}{2}v)=\frac{1}{2}$, for all $x\in{\Z^d}$ and $v\in{V^0}$. We attach to each edge $e=\{x,y\}$ the following interval of length $\frac{1}{2}$: 
\begin{equation}
\label{eq:I_e} 
I_e\stackrel{\text{def.}}{=}\Big\{x+tv_{(x,y)};\, t\in{\Big(0,\frac{1}{2}\Big)}\Big\}=\Big\{y+tv_{(y,x)};\, t\in{\Big(0,\frac{1}{2}\Big)}\Big\},
\end{equation}
which is homeomorphic to an open interval of $\R$ of length $\frac{1}{2},$ and we write $\overline{I_e}=I_e\cup\{x,y\}.$ The cable system $\tilde{\Z}^d$ is then defined by glueing these intervals through their endpoints. $\tilde{\Z}^d\setminus\Z^d$ is now the union of such $I_e$, one for every edge $e\in{E}.$ We extend the definition of the distance $d$ to $\tilde{\Z}^d$ by setting $d(x+tv,x+t'v)=t-t'$ for all $x\in{\Z^d},$ $v\in{V^0}$ and $0\leqslant t'\leqslant t\leqslant\frac{1}{2},$ and for all $x_1\sim x_2$ and $y_1\sim y_2$ in ${\Z^d},$ $z\in{I_{\{x_1,x_2\}}}$ and $z'\in{I_{\{y_1,y_2\}}},$
\begin{equation*}
    d(z,z')=\min_{i,j\in{\{1,2\}}} \big\{ d(z,x_i)+d(x_i,y_j)+d(y_j,z') \big\}.
\end{equation*}
For all $e\in{E}$ and $z_1,\,z_2\in{\overline{I_e}}$ we define $(z_1,z_2)\subset{\tilde{\Z}^d}$ as the open interval in $I_e$ between $z_1$ and $z_2.$ We also define the distance between two subsets $A_1$ and $A_2$ of $\tilde{\Z}^d$ by $d(A_1,A_2)=\inf_{x\in{A_1},y\in{A_2}}d(x,y).$ For $R_1<R_2,$ we introduce the boxes $[R_1,R_2)^d=\{z\in{\tilde{\Z}^d};\, z\in{\overline{I}_{\{x,y\}}},\ \text{with } x_i,y_i\in{[R_1,R_2)} \text{ for all } i =1,\dots, d\}$. The set $\Z^d$ will henceforth be considered as a subset of $\tilde{\Z}^d$ and we will call vertices the elements of $\Z^d.$  

One can define a continuous diffusion $\tilde{X}$ on the cable system $\tilde{\Z}^d,$ via probabilities $\tilde{P}_z,\ z\in{\tilde{\Z}^d},$ with continuous local times with respect to the Lebesgue measure on $\tilde{\Z}^d.$ We now describe this construction from a simple random walk on $\Z^d$ with the help of the excursion process of Brownian motion as in Section 2 of \cite{MR3502602}, and refer to \cite{MR1838739} or Section \nolinebreak2 of \cite{MR3152724} for precise definitions. Let $n$ be the intensity measure of Brownian excursions, see Chapter XII \S 2 in \cite{MR1725357}, and $\lambda_+$ be the Lebesgue measure on $[0,\infty).$ For all $x\in{\Z^d},$ we define under $\tilde{P}_x$  a Poisson point process $\textbf{e}=\sum_{n\in{\N}}\delta_{(e_n,t_n)}$ with intensity measure $n\otimes\lambda_+,$ $(V_n)_{n\in{\N}}$ an i.i.d.\ sequence of uniform variables on $V^0$ independent of $\textbf{e}$ (here and in the sequel $\mathbb{N}=\{ 0,1,2,\dots\}$), and $(Z_n)_{n\in{\N}}$ an independent simple random walk on $\Z^d$ with $Z_0=x$. For any trajectory $e$ in the space of excursions, let $R(e)=\inf\{t>0:\ e(t)=0\}$ be the length of $e,$ and we define for all $n\in{\N}$
\begin{equation*}
    \tau_n \equiv \tau_n(\textbf{e}) :=\sum_{\substack{p\in{\N}\\t_p\leqslant t_n}}R(e_p),\qquad\qquad \tau_n^- \equiv \tau_n^-(\textbf{e}) :=\sum_{\substack{p\in{\N}\\t_p<t_n}}R(e_p), \quad   \text{if } \textbf{e}=\sum_{n\in{\N}}\delta_{(e_n,t_n)},
\end{equation*}
and $T\in{[0,\infty)}$ such that there exists $N\in{\N}$ with $|e_N(T-\tau_N^-)|=\frac{1}{2}$ and for all $p\in{\N}$ such that $t_p<t_N,$ $\sup_{s>0}e_p(s)<\frac{1}{2}.$ In words, $T$ is the first time that the graph obtained by concatenating the excursions in the support of $\textbf{e}$ according to their label $t_n$ reaches height $1/2$ in absolute value.
 For each $x\in{\Z^d},$ we then define under $\tilde{P}_x$ for all $t<\tau_N^-,$
\begin{equation*}
    \tilde{X}_t=x+|e_{n}(t-\tau_{n}^-)|V_n\text{ whenever } \tau_{n}^-\leqslant t\leqslant \tau_{n}.
\end{equation*}
and for all $t\in{[\tau_N^-,T]},$
\begin{equation*}
    \tilde{X}_t=x+|e_{N}(t-\tau_{N}^-)|v_{(Z_0,Z_1)}.
\end{equation*}
Note that, under $\tilde{P}_x,$ $(\tilde{X}_t)_{t\leqslant T}$ is a continuous process on $\bigcup_{y\sim x}\overline{I}_{\{x,y\}},$ and that $\tilde{X}_T=Z_1,$ and we repeat this process after time $T$ starting in $Z_1$ in such a way that, conditionally on $(\tilde{X}_t)_{t\leqslant T},$ the law of $(\tilde{X}_t)_{t\geqslant T}$ is $\tilde{P}_x$-a.s.\ the same as the law of $(\tilde{X}_t)_{t\geqslant0}$ under $\tilde{P}_{Z_1},$ and that the projection of the trajectory of $\tilde{X}$ on $\Z^d$ is $(Z_n)_{n\in{\N}}.$ On an edge, the process $X$ behaves like a Brownian motion, see Chapter XII, Proposition 2.5 in \cite{MR1725357} for a similar construction of the Brownian motion on $\R$ from the Poisson point process of excursions. Finally, for all $x\sim y\in{\Z^d}$ and $z\in{I_{\{x,y\}}},$ we construct $(\tilde{X}_t)_{t\geqslant0}$ under $\tilde{P}_z$ as a Brownian motion beginning in $z$ on $I_{\{x,y\}}$ until either $x$ or $y$ is reached, and then we continue with the previous construction beginning at this vertex.

Under $\tilde{P}_x$ for $x\in{\Z^d},$ the local time in $x$ of $\tilde{X}$ at time $T$ relative to the Lebesgue measure on $\tilde{\Z}^d$ has the same law upon renormalisation as the local time in $0$ of a Brownian motion at the moment it leaves $(-\frac{1}{2},\frac{1}{2}),$ and is thus an exponential variable, see for example Chapter VI, Proposition 4.6 in \cite{MR1725357} for a similar result, with parameter $1,$ see Section 2 of \cite{MR3502602} for details. For all $t\in{[0,\infty]},$ let us denote by $(L^y_t)_{y\in{\tilde{\Z}^d}}$ the local times relative to the Lebesgue measure on $\tilde{\Z}^d$ of $\tilde{X}$ at time $t,$ see Section 2 in \cite{MR3502602}, then for all $x\in{\Z^d},$ $(L^y_{\infty})_{y\in{\Z^d}}$ has the same law under $\tilde{P}_x$ as the field of occupation times of the jump process $X$ on $\Z^d$ under $P_x$ (cf. below \eqref{green}). In particular, we can define for all $x,y\in{\tilde{\Z}^d}$ the Green function 
\begin{equation}
\label{eq:green_cont}
g(x,y)=\tilde{E}_x[L^y_{\infty}],
\end{equation}
and its restriction to $\Z^d$ is the same as the Green function on $\Z^d$ defined in (\ref{green}), so the identical notation does not bear any risk of confusion.

We endow the canonical space $\Omega_0:=C(\tilde{\Z}^d,\R)$ of continuous real-valued functions on $\tilde{\Z}^d$ with the canonical $\sigma$-algebra generated by the coordinate functions $\tilde{\Phi}_x,$ $x\in{\tilde{\Z}^d},$ and let $\tilde{\P}^G$ be the probability on $\Omega_0$ such that, under $\tilde{\P}^G,$

\begin{equation}
\label{tildephi}
\begin{split}
    &(\tilde{\Phi}_x)_{x\in{\tilde{\Z}^d}}\text{ is a centered Gaussian field with}\\
    &\text{covariance function } \tilde{\E}^G[\tilde{\Phi}_x\tilde{\Phi}_y]=g(x,y)\text{ for all }x,y\in{\tilde{\Z}^d},
    \end{split}
\end{equation}
with $g(\cdot, \cdot)$ given by \eqref{eq:green_cont}. With a slight abuse of notation, any random variable $\tilde{\phi}=(\tilde{\phi}_x)_{x\in{\tilde{\Z}^d}}$ on $C(\tilde{\Z}^d,\R)$ with law $\tilde{\P}^G$ under some $\tilde{\P}$ will be called a Gaussian free field on the cable system $\tilde{\Z}^d,$ and it is plain that the restriction of a Gaussian free field on the cable system to $\Z^d$ is a Gaussian free field on $\Z^d,$ so we will often identify $\tilde{\Phi}_x$ with $\Phi_x$ for $x\in{\Z^d}.$ 

Let us recall the simple Markov property for $\tilde{\phi}.$ Let $K\subset\tilde{\Z}^d$ be a compact subset with finitely many components, and let $U=\tilde{\Z}^d\setminus K.$ For all $x\in{\tilde{\Z}^d},$ we define
\begin{equation}
\label{phi^u}
    \tilde{\beta}_x^U =\tilde{E}_x\left[\tilde{\phi}_{\tilde{X}_{T_U}}1_{\{T_U<\infty\}}\right]\quad \text{and} \quad \tilde{\phi}^U_x=\tilde{\phi}_x-\tilde{\beta}_x^U,
\end{equation}
where $T_U:=\inf\{t\geqslant0;\, \tilde{X}_t\notin{U}\},$ with the convention $\inf\emptyset=\infty,$ is the exit time from $U$ of the diffusion $\tilde{X}$ on $\tilde{\Z}^d$. Moreover, for all $x,y\in{\tilde{\Z}^d},$ we define similarly as in \eqref{eq:green_cont} the Green function $g_U(x,y)=\E_x[L^y_{T_U}]$ of the diffusion $\tilde{X}$ under $\tilde{P}_x$ killed when exiting $U.$ Then,
\begin{equation}
\label{markov}
\begin{split}
&(\tilde{\phi}_x^U)_{x\in{\tilde{\Z}^d}}\text{ is a centered Gaussian field with}\\
&\text{covariance function }g_U(x,y)\text{ for all } x,y\in{\tilde{\Z}^d}.
\end{split}
\end{equation}
Furthermore, this field is continuous, vanishes on $K$ and is independent of $\sigma(\tilde{\phi}_z,\ z\in{K}).$ A strong Markov property is also known to hold, but we will not need it here, see Section \nolinebreak1 of \cite{MR3492939} for more details.

Following standard notation, we say that $(B_t)_{t\in{[0,l]}}$ is a \textit{Brownian bridge of length $l>0$ between $x$ and $y$ of a Brownian motion with variance $\sigma^2$ at time $1$} under a probability $\P^B$ if the process
\begin{equation}
W_t:=B_t-\frac{t}{l}y-\Big(1-\frac{t}{l}\Big)x, \quad t \in [0,l],
\end{equation}
is a centered Gaussian field with covariance function
\begin{equation}
\label{bb}
    \E^B\left[W_{s_1}W_{s_2}\right]= \frac{\sigma^2s_1(l-s_2)}{l}\text{ for all }s_1,s_2\in{\left[0,l\right]}\text{ with }s_1\leqslant s_2
\end{equation}
(the process $(W_{l t}/ \sqrt {l \sigma^2})_{t \in [0,1]}$ is a standard Brownian bridge). Let $e\in{E},$ $z_1\neq z_2\in{\overline{I_e}},$ $v\in{V^0}$ and $t\in{(0,\frac{1}{2}]}$ such that $z_2=z_1+tv,$ and let $s_1,s_2\in{[0,t]}$ such that $s_1\leqslant s_2.$ Under $\tilde{P}_{z_1+s_1v},$ until time $T_{(z_1,z_2)},$ the diffusion $\tilde{X}$ behaves like a Brownian motion on $I_e$ beginning at $z_1+s_1v$ until the hitting time of $(z_1,z_2)^c.$ Using Chapter II.11 in \cite{MR1912205} with $s(x)=x,$ and noting that the function $G_0$ defined therein is $\frac{1}{2}g_{(z_1,z_2)},$ we have
\begin{equation*}
    g_{(z_1,z_2)}(z_1+s_1v,z_1+s_2v)=\frac{2s_1(t-s_2)}{t}.
\end{equation*}
The Markov property for the Gaussian free field implies that, under $\tilde{\P}$ (under which $\widetilde{\varphi}$ is a Gaussian free field),
\begin{equation}
\label{bridgeinsideedge}
\Big(\tilde{\phi}_{z_1+sv}-\frac{t-s}{t}\tp_{z_1}-\frac{s}{t}\tp_{z_2}\Big)_{s\in{[0,t]}}
\end{equation}
is a centered Gaussian field with covariance function $(g_{(z_1,z_2)}(z_1+s_1v,z_2+s_2v))_{s_1,s_2\in{[0,t]}},$ and is independent of $\sigma(\tilde{\phi}_z,\,z\in{\tilde{\Z}^d\setminus(z_1,z_2))}.$ Thus, it is a Brownian bridge of length $t$ between $0$ and $0$ of a Brownian motion with variance $2$ at time $1.$ In particular, knowing $\tilde{\phi} \restriction \Z^d,$ the Gaussian free field on the edges $((\tilde{\phi}_z)_{z\in{I_e}})_{e\in{E}}$ is an independent family of random processes such that, for each $x\sim y\in{\Z^d},$ the process $(\tilde{\phi}_z)_{z\in{I_{\{x,y\}}}}$ has the same law as a Brownian bridge of length $\frac{1}{2}$ between $\tilde{\phi}_{x}$ and $\tilde{\phi}_{y}$ of a Brownian motion with variance 2 at time 1, as mentioned in Section 2 of \cite{MR3502602} or in Section 2.2 of \cite{LupuSabotTarres}. More precisely, let 
\begin{equation}
\label{brownianbridgeontheedges}
    B_t^e=\tilde{\phi}_{x_e+tv_{(x_e,y_e)}}-2t\tp_{y_e}-(1-2t)\tp_{x_e},\text{ for all }t\in{\left[0, 1/2\right]}\text{ and }e\in{E},
\end{equation}
where we have given an (arbitrary) orientation $(x_e,y_e)$ for each edge $e=\{x_e,y_e\}\in{E}.$ Then, under $\tilde{\P},$ $(B^e)_{e\in{E}}$ is a family of independent Brownian bridges of length $\frac{1}{2}$ between 0 and 0 of a Brownian motion with variance 2 at time 1. Note that this provides an explicit (and simple) construction of a Gaussian free field on the cable system starting from the Gaussian free field $(\phi_x)_{x\in{\Z^d}}$ on $\Z^d$: if one links independently each $x\sim y\in{\Z^d}$ via a Brownian bridge on $I_{\{x,y\}}$ of length $\frac{1}{2}$ between $\phi_x$ and $\phi_y$ of a Brownian motion with variance $2$ at time $1,$ then the resulting process is a Gaussian free field on the cable system. In view of this construction, we will later need the following result on the probability that the maximum of a Brownian bridge exceeds some value $M$ (see e.g.\ \cite{MR1912205}, Chapter IV.26).
\begin{Lemme}
\label{Brownianbridge}
Let $x,\,y$ be two real numbers, $M\geqslant \max(x,y)$ and, under $\P^B,$ $(B_t)_{t\in{[0,l]}}$ a Brownian bridge of length l between $x$ and $y$ of a Brownian motion with variance $\sigma^2$ at time 1. One has
\begin{equation}
\label{sup}
    \P^B\Big(\sup_{t\in{[0,l]}}B_t> M\Big)=\exp\left(-\frac{2(M-x)(M-y)}{l\sigma^2}\right).
\end{equation}
\end{Lemme}
Let us now turn to the definition of random interlacements on $\tilde{\Z}^d,$ as in \cite{MR3502602} or \cite{MR3492939}. The usual definition of random interlacements on $\Z^d$, see, for example, \cite{MR2680403} or the monograph \cite{MR3308116}, can be adapted to define a Poisson point process $\tilde{\omega}$ on $\tilde{W}^*\times[0,\infty),$ where $\tilde{W}^*$ is the space of doubly infinite trajectories on $\tilde{\Z}^d$ modulo time-shift, endowed with its canonical $\sigma$-algebra, and where $[0,\infty)$ describes labels of the trajectories. Recall the law $(\tilde{P}_z)_{z\in{\tilde{\Z}^d}}$ of the diffusion $\tilde{X}$ on the cable system, started at $z\in{\tilde{\Z}^d}$. The intensity measure of  $\tilde{\omega}$ is characterized as follows: for some $N_1,N_2\in{\Z}$ with $N_1\leqslant N_2,$ let $\tilde{K}:=[N_1,N_2]^d\cap\tilde{\Z}^d,$ let $K:=\tilde{K}\cap\Z^d,$ let $\tilde{\omega}^u$ be the point process which consists of the trajectories in $\tilde{\omega}$ with label at most $u>0$, and let $\tilde{\omega}^u_{\tilde{K}}$ be the point process comprising the forward trajectories of $\tilde{\omega}^u$ hitting $\tilde{K}$ and beginning at the first time $\tilde{K}$ is reached. Then $\tilde{\omega}^u_{\tilde{K}}$ is a Poisson point process with intensity measure $u\tilde{P}_{e_{K}}= u \sum_{x \in K} e_K(x)\tilde{P}_x$, where $e_{K}$ is the usual equilibrium measure of $K$ on $\Z^d,$ as mentioned in \cite{MR3492939}. One can also construct the random interlacement process $\tilde{\omega}^u$ at level $u>0$ on the cable system from the corresponding interlacement process $\omega^u$ on $\Z^d$ by adding independent Brownian excursions on the edges for every trajectory in the support of $\omega^u$ in the same fashion as one can construct the diffusion $\tilde{X}$ from a simple random walk on $\Z^d.$

One then defines
\begin{equation}
\label{tildel}
( \tilde{\ell}_{z,u})_{z\in{\tilde{\Z}^d}}\text{ the field of local times of random interlacements, for $u>0$,}
\end{equation}
as the sum of the local times of each of the trajectories in the support of $\tilde{\omega}^u.$ The restriction of these local times to $\Z^d$ has the same law as the occupation times $(\ell_{x,u})_{x \in \Z^d}$ for random interlacements on $\Z^d$ alluded to in the introduction, cf.\ above \eqref{Iu}. The random interlacement set is defined as
\begin{equation}
\label{itildeu}
    \tilde{\I}^u=\{x\in{\tilde{\Z}^d};\,  \tilde{\ell}_{x,u}>0\},
\end{equation} 
which is an open connected subset of $\tilde{\Z}^d$. Note that $\{x\in{\Z^d};\, x\in{\tilde{\I}^u}\}$ has the same law as $\I^u,$ cf.\ \eqref{Iu}. 

We also recall the following formula for the Laplace transform of $(\ell_{x,u})_{x \in \Z^d}$, see for instance \cite{MR2892408}, (1.9)--(1.11) or Remark 2.4.4 in \cite{MR3050507}: for all $V: \Z^d \to \mathbb{R}$ with finite support $K \subset \mathbb{Z}^d$ and satisfying
\begin{equation}
\label{eq:lap1}
\Vert GV \Vert_{\infty} < 1,  \text{ where } (GV)f(x) = \sum_{y \in \Z^d}g(x,y)V(y)f(y)\text{ for all }f\in{\ell^{\infty}(\Z^d)},
\end{equation}
with $g(\cdot, \cdot)$ as in \eqref{green}, and where $\Vert \cdot \Vert_{\infty}$ denotes the operator norm on $\ell^{\infty}(\Z^d)\rightarrow \ell^{\infty}(\Z^d),$ one has
\begin{equation}
\label{eq:lap2}
\begin{split}
\tilde{\E}^I \Big[ \exp \Big\{ \sum_{x\in \Z^d} V(x) \tilde{\ell}_{x,u} \Big\} \Big] 
&\ \ = \exp\{ u \langle V, (I- GV)^{-1} 1 \rangle_{L^2(\mathbb{Z}^d)}\}\\
&\Big(= \exp \Big\{ u \sum_{x \in \Z^d} V(x)  \sum_{ n \geqslant 0} (GV)^n1(x) \Big\} \Big).
\end{split}
\end{equation}

Random interlacements are useful in the study of the Gaussian free field on the cable system $\tilde{\Z}^d$ because of the existence of a Ray-Knight-type isomorphism theorem proved in Proposition 6.3 of \cite{MR3502602}, see also (1.30) in \cite{MR3492939}.  

\begin{The}

\label{iso}
For each $u>0,$ there exists a coupling $\tilde{\P}^u$ between two Gaussian free fields $\tilde{\phi}$ and $\tilde{\gamma}$ and a random interlacement process $\tilde{\omega}$ on the cable system $\tilde{\Z}^d$ (i.e., under $\tilde{\P}^u,$ the law of $\tilde{\phi}$ and $\tilde{\gamma}$ is $\tilde{\P}^G$ each, and the law of $\tilde{\omega}$ is the same as under $\tilde{\P}^I$) such that $\tilde{\gamma}$ and $\tilde{\omega}$ are independent, and $\tilde{\P}^u$-a.s.,
\begin{equation}
\label{isoeq}
   \frac{1}{2}\left(\tilde{\phi}_x+\sqrt{2u}\right)^2= \tilde{\ell}_{x,u}+\frac{1}{2}\tilde{\gamma}_x^2,\quad\text{ for all }x\in{\tilde{\Z}^d},
\end{equation}
where $( \tilde{\ell}_{x,u})_{x\in{\tilde{\Z}^d}}$ is the field of local times of the random interlacements process $\tilde{\omega}$ at level~$u,$ cf.\ \eqref{tildel}.
\end{The}

This coupling will be essential for the proof of Theorem \ref{mainresult}. In particular, we are going to use results from the theory of random interlacements, along with the coupling (\ref{isoeq}), to deduce certain properties of the level sets of the Gaussian free field. For now, let us note that $\tilde{\P}^u$-a.s.\ on $\tilde{\I}^u,$ cf.\ \eqref{itildeu}, one has $|\tilde{\phi}+\sqrt{2u}|>0$, where $\tilde{\phi}$ refers to the Gaussian free field from the coupling in Theorem \ref{iso}. Since $\tilde{\I}^u$  is connected (and unbounded, by construction) and since $x\in{\tilde{\Z}^d}\mapsto\tilde{\phi}_x$ is continuous, either $\tilde{\phi}_x>-\sqrt{2u}$ for all $x\in{\tilde{\I}^u},$ or $\tilde{\phi}_x<-\sqrt{2u}$ for all $x\in{\tilde{\I}^u}.$ But Proposition 5.5 in \cite{MR3502602}, cf.\ also \eqref{eq:1.20} above, implies that the set $\{x\in{\tilde{\Z}^d}; \, \tilde{\phi}_x<0\}$, which contains $\{x\in{\tilde{\Z}^d}; \, \tilde{\phi}_x<-\sqrt{2u}\}$, only has bounded components, hence
\begin{equation}
\label{htilde=0}
\tilde{\P}^u-\text{a.s.,}\ \forall x\in{\tilde{\I}^u},\ \tilde{\phi}_x>-\sqrt{2u}. 
\end{equation}
In particular, this means that the negative (upper) level sets percolate on $\tilde{\Z}^d$, see~\eqref{unboun}.

\section{Connectivity and a large deviation inequality for \texorpdfstring{$\tilde{\mathcal{I}}^u$}{random interlacements}
}
\label{strong}

The following result, which is proved over the next two sections, is essentially a refinement of \eqref{htilde=0}, which allows us to truncate $\widetilde{\Phi}$, cf. \eqref{tildephi}, at sufficiently large heights. This important technical step will be helpful in dealing with the fact that $\widetilde{\Phi}$ is a priori unbounded on sets of interest.
\begin{The}
\label{percphitilde}
For each $h_0>0$, there exist positive constants $C_0$ and $c_0,$ only depending on $d$ and $h_0,$ with $C_0h_0^{-c_0}\geqslant1$ such that, for all $0 < h\leqslant h_0,$ with 
\begin{equation}
\label{K}
K(h)=\sqrt{\log\left(\frac{C_0}{h^{c_0}}\right)},
\end{equation}
there exists $L_0 = L_0(h)>0$ such that $\tilde{\P}^G$-a.s.\ the set 
\begin{equation}
\label{A''h}
    \tilde{A}_h(\tilde\Phi) \stackrel{\text{def.}}{=}\left\{x\in{\tilde{\Z}^d\setminus{\Z^d}};\, \tilde{\Phi}_{x}\geqslant -h\right\}\cup\Big\{x\in{\Z^d};\, \forall\,v\in{V^0},\, \forall t\in{\Big[0,\frac{1}{2}\Big]},\, |\tilde{\Phi}_{x+tv}|\leqslant K(h)\Big\}
\end{equation}
contains an unbounded connected component in the thick slab $\tilde{\Z}^2\times[0,2L_0)^{d-2}.$
\end{The}

Note that, since $\tilde{\Phi}$ is continuous, asserting that $\tilde{A}_h(\tilde{\Phi})$ has an unbounded component is tantamount to saying that there exists an infinite path in the set $\{x\in{\tilde{\Z}^d};\, \tilde{\Phi}_{x}\geqslant -h\},$ and that in addition, for every $y\in{\tilde{\Z}^d}$ at distance less than $\frac{1}{2}$ from a vertex on this path, $|\tilde{\Phi}_y|\leqslant K(h)$ holds. In particular, the set    $\{x\in{\tilde{\Z}^d};\,-h\leqslant\tilde{\Phi}_x\leqslant K(h)\}$ also contains an unbounded connected component in the thick slab $\tilde{\Z}^2\times[0,2L_0)^{d-2},$ as stated in \eqref{1:30}. In order to be able to prove Theorem \ref{mainresult} with the help of Theorem \ref{percphitilde} in Section \ref{proofofthemainresult}, the key property of $K(h)$ in \eqref{K} is that 
\begin{equation}
\label{Kcondition}
hK(h)\rightarrow0,  \text{ as } h\searrow0,
\end{equation}
 see in particular the proof of Lemma \ref{comparingEandF}. The proof of Theorem \ref{percphitilde} will involve an application of the isomorphism \eqref{isoeq}, and therefore hinges on a corresponding statement ``in the world of random interlacements,'' see Proposition \ref{truncated negative level sets} at the beginning of the next section. The proof of the latter requires some preliminary results on the geometry of $\tilde{\I}^u,$ which we gather now. The dependence of these results on $u$ needs to be precise enough to later deduce \eqref{Kcondition} when transferring Proposition \ref{truncated negative level sets} back to the Gaussian free field.

In the remainder of this section, we consider, under $\tilde{\P}^I$, and for each $u>0$, the random interlacement set $\tilde{\I}^u$ at level $u$ on the cable system, see \eqref{itildeu}, and $( \tilde{\ell}_{x,u})_{x\in{\tilde{\Z}^d}}$ the field of local times of the underlying interlacement process $\tilde\omega^u$, see \eqref{tildel}. The following lemma asserts that $\tilde{\mathcal{I}}^u$ is typically well-connected.

\begin{Lemme}
\label{connectivity}
Let $d\geqslant3,$ $\epsilon\in{(0,1)}$ and  $u_0>0.$ There exist constants $c=c(d,\epsilon,u_0)$ and $C=C(d,\epsilon,u_0)$ such that for all $u\in{(0,u_0]}$ and $R\geqslant1,$
\begin{equation}
\label{eq:connectivity}
\tilde{\P}^I\bigg(\bigcap_{x,y\in\tilde{\I}^u\cap[0,R)^d}\{x\leftrightarrow y\ in\ \tilde{\I}^u\cap[-\epsilon R,(1+\epsilon)R)^d\}\bigg) 
\geqslant 1-C\exp\left(-cR^{1/7}u\right),
\end{equation}
where, for measurable $A\subset\tilde{\Z}^d,$ the event $\{x\leftrightarrow y\ in\ \tilde{\I}^u\cap A\}$ refers to the existence of a continuous path in the subset $\tilde{\I}^u\cap A$ of the cable system connecting $x$ and $y$. 
\end{Lemme}

This property is essentially known, see for instance Proposition 1 of \cite{MR2819660} or Lemma 3.1 in \cite{MR3024098}. However, we need to keep careful track of the dependence of error terms on the intensity $u$. 
For the reader's convenience, we have included a proof of Lemma \ref{connectivity} in the Appendix. 

Next, we will need to know how 
much time the trajectories of random interlacements typically spend
 in a large box with sufficiently high precision. This can be conveniently formulated in terms of a large deviation inequality for the local times. 

\begin{Lemme}
\label{largedeviationlemma}
Let $d\geqslant3$ and $\epsilon\in{(0,1)}$. There exist constants $c=c(d,\epsilon)$ and $C=C(d,\epsilon)$ such that for all $u>0$ and $R\geqslant1$,
\begin{equation}
\label{1+eps}
    \tilde{\P}^I\bigg( \, \Big|\frac{1}{R^d}\sum_{x\in{[0,R)^d\cap\Z^d}}{ \tilde{\ell}_{x,u}} - u \Big| > \varepsilon \cdot u \bigg)
    \leqslant C\exp\left(-cR^{d-2}u\right).
\end{equation}
\end{Lemme}
\begin{proof} Abbreviate $B_R = [0,R)^d\cap\Z^d$ and, for $\lambda > 0$, let $V(x)=  (|B_R|)^{-1} \lambda 1_{\{ x\in B_R \}}$. It follows, cf.\ \eqref{eq:lap1} for notation, that there exists $K_1<\infty$ such that for all $f\in{\ell^{\infty}(\Z^d)}$ and $x \in \Z^d$,
\begin{equation}
\label{eq:lap3}
|(GV)f(x)| =  \lambda \Big| \sum_{y \in B_R} \frac{g(x,y)f(y)}{|B_R|} \Big| \leq K_1 \lambda R^{2-d}\Vert f\Vert_{\ell^{\infty}(\Z^d)},
\end{equation}
using that $g(x,y) \leq C' |x-y|^{2-d}$, for $x,y \in \Z^d$. Hence, for $\lambda= \lambda_0 R^{d-2}$, with $\lambda_0 < K_1^{-1}$, one obtains that 
$\Vert GV\Vert_{\infty}< 1$, for all $R \geq 1$. In view of \eqref{eq:lap3}, applying Markov's inequality and using \eqref{eq:lap2} then yields, for all $\lambda_0 < K_1^{-1}$ and $R \geq 1$,
\begin{align*}
    \tilde{\P}^I\Big(\,\frac{1}{R^d}\sum_{x\in{[0,R)^d\cap\Z^d}}{ \tilde{\ell}_{x,u}}>(1+\epsilon) u\Big)
    &\leqslant \exp \Big\{ - \lambda_0 R^{d-2} u \,  \Big( \big(1+\varepsilon\big) -  \big( 1 + \sum_{ n \geqslant 1} (K_1 \lambda_0)^n \big ) \Big) \Big\}.
\end{align*}
The right-hand side is bounded from above by $C\exp\left(-cR^{d-2}u\right)$ upon choosing $\lambda_0 (\varepsilon) < K_1^{-1}$ small enough such that $\sum_{ n \geqslant 1} (K_1 \lambda_0)^n \leq \varepsilon/2$. In a similar fashion one bounds for $V$, $\lambda$ as above, 
\begin{align*}
    \tilde{\P}^I\Big(\,\frac{1}{R^d}\sum_{x\in{[0,R)^d\cap\Z^d}}{ \tilde{\ell}_{x,u}}<(1 - \epsilon) u\Big)
    &=  \tilde{\P}^I\Big(\,\exp \Big\{  - \sum_{x\in \Z^d} V(x) \tilde{\ell}_{x,u} \Big\}  >  e^{-(1-\varepsilon)\lambda u} \Big) \\ 
   & \leqslant  \exp \Big\{ \lambda_0 R^{d-2} u \,  \Big( \big(1-\varepsilon\big) -   \big(1-\sum_{n\geqslant1}{(K_1\lambda_0)^n}\big) \Big) \Big\},
\end{align*}
from which \eqref{1+eps} readily follows.
\end{proof}

As a direct application of Lemmas \ref{connectivity} and \ref{largedeviationlemma}, we derive lower bounds for the probabilities of the following events. 

\begin{Def}
\label{defEandF}
For all $u,\, u'>0,$ and integer $R\geq 1$, the events $E_{R}^{u,u'}$ and $F_{R}^{u,u'}$ are defined as follows:
\begin{enumerate}[(a)]
\item $E_{R}^{u,u'}$ occurs if and only if for each $e\in\{0,1\}^d,$ the set $(eR+[0,R)^d)\cap\tilde{\I}^u$ contains a connected component $\mathcal{A}_{e}$ such that
        \begin{equation*}
        \sum_{y\in{\mathcal{A}_{e}\cap\Z^d}}{ \tilde{\ell}_{y,u}}>\frac{3}{4}u'R^d,
        \end{equation*}
 and such that the components $(\mathcal{A}_{e})_{e\in\{0,1\}^d}$ are all connected in $\tilde{\I}^u\cap [0,2R)^d.$
\item $F_{R}^{u,u'}$ occurs if and only if for all $e\in\{0,1\}^d,$
\begin{equation*}
        \sum_{y\in{(eR+[0,R)^d)\cap\Z^d}}{ \tilde{\ell}_{y,u}}<\frac{5}{4}u'R^d.
\end{equation*}    
\end{enumerate}
\end{Def}

Note that for fixed $u'>0,$ and positive integer $R,$ the events $(E_R^{u,u'})_{u>0}$ are increasing,  i.e., there exists a measurable and increasing function $f_R^{u'}:[0,\infty)^{\tilde{\Z}^d}\rightarrow\{0,1\}$ such that $1_{E_R^{u,u'}}=f_R^{u'}(\tilde{\ell}_{.,u})$ for all $u>0,$ and that the events $(F_R^{u,u'})_{u>0}$ are decreasing, i.e., the events $((F_R^{u,u'})^c)_{u>0}$ are increasing. The following consequence of Lemmas \ref{connectivity} and \ref{largedeviationlemma} is tailored to our purposes in the next section. 

\begin{Cor}
\label{E0uandF0u}
Let $d\geqslant3$ and $u_0>0.$ There exist $\delta\in{(0,1)},$ positive and finite constants $C=C(d,u_0)$ and $c=c(d,u_0)$ such that for all $u\in{(0,u_0]}$ and $R\geqslant1,$ 
\begin{equation}
\label{E0u}
\tilde{\P}^I\big(E_R^{u(1-\delta),u}\big)\geqslant 1-C\exp\left(-cR^{1/7}u\right)
\end{equation}
and 
\begin{equation}
\label{F0u}
\tilde{\P}^I\big(F_R^{u(1+\delta),u}\big)\geqslant 1-C\exp\left(-cR^{d-2}u\right).
\end{equation}
\end{Cor}

\begin{proof}
Let $\delta=\frac{1}{6}$. We begin with \eqref{F0u}. In view of Definition \ref{defEandF}, it follows from Lemma~\ref{largedeviationlemma} applied with $u(1+\delta)$ instead of $u$ and translation invariance that for all $e\in{\{0,1\}^d}$ and $u > 0$,
\begin{align*}
   & \tilde{\P}^I\left(\, \sum_{x\in{(eR+[0,R)^d)\cap\Z^d}}{ \tilde{\ell}_{x,u(1+\delta)}}<\frac{5}{4}R^d u\right) \\
    &\qquad \geqslant\tilde{\P}^I\left(\frac{1}{R^d}\sum_{x\in{(eR+[0,R)^d)\cap\Z^d}}{ \tilde{\ell}_{x,u(1+\delta)}}< \frac{15}{14}u(1+\delta)\right) \stackrel{\eqref{1+eps}}{\geqslant} 1-C\exp(-cR^{d-2}u),
\end{align*}
which is \eqref{F0u}.

In order to obtain \eqref{E0u}, fix any two constants $\epsilon=\epsilon(d)\in{(0,1)}$ and $\mu=\mu(d)\in{(0,1)}$ in such a way that $(1-8\epsilon)^d(1-\mu)(1-\delta)=\frac{3}{4}.$ For all $e\in{\{0,1\}^d},$ we define the inner boxes $\mathtt{B}_e(\eps)=eR+[2\lfloor\epsilon R\rfloor,R-2\lfloor\epsilon R\rfloor)^d$. It is sufficient to prove \eqref{E0u} for $R$ satisfying $\eps R\geqslant1,$ which we now tacitly assume. We then have $|\mathtt{B}_e(\eps)\cap\Z^d| \cdot (1-\mu)(1-\delta)\geqslant\frac{3}{4}R^d$, where $|A|$ denotes the cardinality of $A\subset\Z^d.$ According to Lemma \ref{largedeviationlemma}, 
\begin{equation*}
\begin{split}
    \tilde{\P}^I\left(\sum_{x\in{\mathtt{B}_e(\eps)\cap\Z^d}}{ \tilde{\ell}_{x,u(1-\delta)}}>\frac{3}{4}R^d u\right)
    & \geqslant\tilde{\P}^I\left(\frac{1}{|\mathtt{B}_e(\eps)\cap\Z^d|}\sum_{x\in{\mathtt{B}_e(\eps)\cap\Z^d}}{ \tilde{\ell}_{x,u(1-\delta)}}> (1-\mu)u(1-\delta)\right) \\
    &\geqslant 1-C\exp\left(-cR^{d-2}u\right).
\end{split}
\end{equation*}
We now define $\mathcal{A}_{e}^1=\mathtt{B}_e(\eps)\cap\tilde{\I}^{u(1-\delta)}.$ According to Lemma \ref{connectivity}, for every $e\in{\{0,1\}^d},$ all the vertices of $\mathcal{A}_{e}^1$ are connected in $\tilde{\I}^{u(1-\delta)}\cap (eR+[\lfloor\epsilon R\rfloor,R-\lfloor\epsilon R\rfloor)^d)$ with probability at least $1-C\exp(-cR^{1/7}u),$ and on the corresponding event we define $\mathcal{A}_{e}\subset \tilde{\I}^{u(1-\delta)}\cap(eR+[\lfloor\epsilon R\rfloor,R-\lfloor\epsilon R\rfloor)^d)$ such that $\mathcal{A}_{e}^1\subset \mathcal{A}_{e}$ and $\mathcal{A}_{e}$ is connected. 

Still according to Lemma \ref{connectivity}, all the $\mathcal{A}_{e}$ for $e\in{\{0,1\}^d}$ are connected with each other in $\tilde{\I}^{u(1-\delta)}\cap [0,2R)^d$ with probability at least $1-C\exp(-cR^{1/7}u),$ which gives (\ref{E0u}).
\end{proof}

Since the events $E_{R}^{u,u'}$ and $F_{R}^{u,u'}$ are defined in terms of local times and not in terms of the occupation field $\big(1_{\{x\in{\tilde{\I}^u}\}}\big)_{x\in{\Z^d}}$, we now give a slightly different version of the decoupling inequality presented in \cite{MR3420516} valid for the local times on the cable system. This inequality will later enable us to use $E_R^{u,u'}$ and $F_R^{u,u'}$ as seed events of a suitable multi-scale argument. In what follows, let $\tilde{\mathbf{Q}}^u$, $u>0$, be the law on $\overline{\Omega}=[0,\infty)^{\tilde{\Z}^d}$ of the local times $( \tilde{\ell}_{x,u})_{x\in{\tilde{\Z}^d}}$ of random interlacements on the cable system $\tilde{\Z}^d$, and let $(p_x)_{x\in{\tilde{\Z}}^d}$ denote the canonical coordinate functions on $\overline{\Omega},$ i.e.,\ for all $f\in{\overline{\Omega}}$ and $x\in{\tilde{\Z}^d},$ $p_x(f)=f(x).$

\begin{The}
\label{couplingforinterlacements} Let $\tilde{A}_1$ and $\tilde{A}_2$ be two measurable non-intersecting subsets of $\tilde{\Z}^d$. Assume that $s:=d(\tilde{A}_1,\tilde{A}_2)\geqslant1,$ and that the minimum $r$ of the diameters of $\tilde{A}_1$ and $\tilde{A}_2$ is finite. Then there exist $\kappa_0(d)$ and $\kappa_1(d)$ such that for all $u>0$ and $\epsilon\in{(0,1)}$,
for any functions $f_i:\overline{\Omega}\rightarrow[0,1]$ which are $\sigma(p_x,\ x\in{\tilde{A}_i})$ measurable for each $i\in{\{1,2\}},$ and which are both increasing or both decreasing,
    \begin{equation}
    \label{decinc}
    \tilde{\mathbf{Q}}^u[f_1f_2]\leqslant \tilde{\mathbf{Q}}^{u(1\pm \eps)}[f_1]\tilde{\mathbf{Q}}^{u(1\pm\eps)}[f_2]+\kappa_0(r+s)^d\exp(-\kappa_1\epsilon^2us^{d-2}),
    \end{equation}
where the plus sign corresponds to the case where the $f_i$'s are increasing and the minus sign to the case where the $f_i$'s are decreasing.
\end{The}
\begin{proof}
Let $A_1$ and $A_2$ be the smallest subsets of $\Z^d$ such that for all $i\in{\{1,2\}},$ and all $x\in{\tilde{A}_i},$ there exist $y,z \in A_i$ such that $x \in \overline{I}_{\{y,z\}}.$ Note that $\tilde{A}_i\cap\Z^d\subset A_i$. Since $d(\tilde{A}_1,\tilde{A}_2)\geqslant1$, the sets $A_1$ and $A_2$ are not intersecting (recall that the distance between two neighbors of $\Z^d$ is $\frac{1}{2}).$ For two measures $\mu_1$ and $\mu_2,$ we say that $\mu_1\leqslant\mu_2$ if $\mu_2-\mu_1$ is a non-negative measure. The proof of the main decoupling result, Theorem 2.1 in \cite{MR3420516}, see in particular Section 5 therein, implies that, for each $u>0,$ there exists a coupling $\Q_u^I$ between the random interlacement process $\omega$ on $\Z^d$ and two independent Poisson point processes $\omega_1$ and $\omega_2$ having the same law as $\omega,$ such that, for $B\subset\Z^d,$ denoting by $(\omega^u)_{|B}$ the point process consisting of the restriction to $B$ of the trajectories in $\omega^u$ which hit $B$,
\begin{equation}
\label{couplinghati}
    \Q_u^I\big[(\omega_i^{u(1-\eps)})_{|A_i}\leqslant(\omega^u)_{|A_i}\leqslant(\omega_i^{u(1+\eps)})_{|A_i},\;i=1,2\big]\geqslant1-\kappa_0(r+s)^d\exp(-\kappa_1\epsilon^2us^{d-2}).
\end{equation}
For each $u>0$ and $i\in{\{1,2\}},$ under an extended probability $\tilde{\Q}_u^I,$ one then constructs an interlacement process $\tilde{\omega}_i^{u(1-\eps)}$ at level $u(1-\eps)$ on the cable system by adding independent Brownian excursions on the edges for every trajectory in the support of the random interlacement process $\omega_i^{u(1-\eps)}$ on $\Z^d,$ as in the construction of the diffusion $\tilde{X},$ see the beginning of Section \ref{cablesystem}.

We now construct a random interlacement process $\tilde{\omega}^u$ at level $u$ using $\tilde{\omega}_i^{u(1-\eps)}$ and $\omega^u.$ Its trajectories are the trajectories of $\tilde{\omega}_i^{u(1-\eps)}$ which have a projection on $\Z^d$ already contained in $\omega^u$ (i.e., all the trajectories of $\tilde{\omega}_i^{u(1-\eps)}$ on the event in \eqref{couplinghati}) and the trajectories of $\omega^u$ which are not already in $\omega_i^{u(1-\eps)}$, lifted to $\tilde \Z^d$ using additional independent (of $\tilde{\omega}_i^{u(1-\eps)}$ and $\omega^u$) Brownian excursions on the edges. We repeat this construction to obtain a random interlacement process $\tilde{\omega}_i^{u(1+\eps)}$ at level $u(1+\eps)$ in a similar way from $\tilde{\omega}^u$ and $\omega_i^{u(1+\eps)}.$ Then, an analogue of \eqref{couplinghati} holds for these processes $\tilde{\omega}^u,$ $\tilde{\omega}^{u(1-\eps)}_i$ and $\tilde{\omega}_i^{u(1+\eps)}$ under $\tilde{\Q}_u^I$. In particular, denoting by $ \tilde{\ell}_{x,u},$ $ \tilde{\ell}^{\,i}_{x,u(1-\eps)}$ and $ \tilde{\ell}^{\,i}_{x,u(1+\eps)}$ their respective local time fields on the cable system, see \eqref{tildel}, it follows that
\begin{equation}
\label{coup100}
    \tilde{\Q}_u^I\big[ \tilde{\ell}^{\,i}_{x,u(1-\eps)}\leqslant \tilde{\ell}_{x,u}\leqslant \tilde{\ell}^{\,i}_{x,u(1+\eps)},\;x\in{\tilde{A}_i},\;i=1,2\big]\geqslant1-\kappa_0(r+s)^d\exp(-\kappa_1\epsilon^2us^{d-2})
\end{equation}
The inequalities in  \eqref{decinc} are a direct consequence of \eqref{coup100}.
\end{proof}

\section{Percolation for the truncated level set}
\label{truncatedGFF}

In this section, we prove Theorem \ref{percphitilde}: for each $h>0,$ there exists a finite constant $K(h)$ such that the level set of the Gaussian free field on the cable system truncated above level $-h$ and below level $K(h)$ contains an unbounded connected component. We will actually show a similar statement for random interlacements and use the coupling from Theorem \ref{iso} to obtain Theorem \ref{percphitilde}. The corresponding statement for random interlacement, see Proposition \ref{truncated negative level sets}, essentially asserts that one can intersect the continuous interlacement set $\tilde{\I}^u$, the set $\{x\in{\Z^d};\, |\phi_x|<K\}$ and a Bernoulli family on the edges with parameter $p$ and still retain an unbounded connected component in $\tilde{\Z}^d$ for sufficiently large $K$ and $p$ close enough to $1$. The proof of this statement bears similarities to the proof of Theorem 2.1 in \cite{MR3024098}, where it is shown that the intersection of $\I^u$ and a Bernoulli family with parameter $p$ on $\Z^d,$ not necessarily independent from $\I^u,$ contains an infinite connected component in $\Z^d$ for large enough $p$.

Henceforth, for a given $p\in{(0,1)}$ (and $d\geqslant 3$), let $\tilde{\Q}^p$ be any coupling between a Gaussian free field $\tilde{\phi},$ a random interlacement process $\tilde{\omega}$ and a family of independent Bernoulli random variables on the edges $\mathcal{B}^p=(\theta_e^{p})_{e\in{E}}$ with parameter $p,$ i.e., 
\begin{equation}
\label{eq:Qdef}
\begin{split}
\text{under $\tilde{\Q}^p,$ } &\text{the law of $(\tilde{\phi}_x)_{x\in{\tilde{\Z}^d}}$ is $\tilde{\P}^G,$ the process $\tilde{\omega}$ has the same law as under $\tilde{\P}^I$}\\
&\text{and $(\theta_e^{p})_{e\in{E}}$ is an i.i.d.\ family of $\{0,1\}$-valued random variables with}\\
&\text{$\tilde{\Q}^p(\theta_e^p=1)=p$ for each $e\in{E}.$} 
\end{split}
\end{equation}
In particular, $\tilde{\phi},$ $\tilde{\omega}$ and $\mathcal{B}^p$ need not be independent, and in fact, we will later use a coupling such that \eqref{htilde=0} holds. For any level $u>0,$ we define the random interlacement set $\tilde{\I}^u$ as in \eqref{itildeu} and the local times $(\tilde{\ell}_{x,u})_{x\in{\tilde{\Z}^d}}$ as in \eqref{tildel} in terms of $\tilde{\omega}$. We further denote by $\phi$ the restriction of $\tilde{\phi}$ to $\Z^d$ and by $\I^u$ the restriction of $\tilde{\I}^u$ to $\Z^d.$ 

\begin{Prop}
\label{truncated negative level sets}
$(d \geqslant 3, \, u_0 > 0,\, \eqref{eq:Qdef} )$

\medskip
\noindent There exist positive constants $C_1,$ $c_1,$ $C'_1$ and $c'_1,$ only depending on $d$ and $u_0,$ satisfying $C_1u_0^{-c_1}\geqslant1$ and $C'_1u_0^{c'_1}<1,$ such that for all $u\in{(0,u_0]},$ with
\begin{equation}
\label{Kpp'}
\tilde{K}(u) \stackrel{\textnormal{def.}}{=}\sqrt{\log\left(\frac{C_1}{u^{c_1}}\right)}\quad \text{ and } \quad p(u) \stackrel{\textnormal{def.}}{=}1-C'_1u^{c'_1},
\end{equation}
there exists $L_0(u)>0$ such that, if $p\in[p(u),1],$ then $\tilde{\Q}^{p}$-a.s.\ the set 
\begin{equation}
\label{A'u}
    \tilde{A}'_{u,p}=\left(\tilde{\I}^u \setminus {\I}^u\right)\cup\left\{x\in{\I^u};\, |\phi_{x}|\leqslant \tilde{K}(u)\text{ and }\forall\, y\sim x,\, |\phi_{y}|\leqslant \tilde{K}(u)\text{ and }\theta^{p}_{\{x,y\}}=1\right\}
\end{equation}
contains an unbounded connected component in the thick slab $\tilde{\Z}^2\times[0,2L_0(u))^{d-2}.$
\end{Prop}
We now comment on \eqref{A'u}. First, note that $\I^u\subset\Z^d,$ so saying that $\tilde{A}'_{u,p}$ contains an unbounded connected component implies that $\tilde{A}'_{u,p}\cap\Z^d$ contains an infinite path such that all the edges of this path are in $\tilde{\I}^u,$ and for all vertices $x$ on this path and all $y\sim x$, $|\phi_x|\leqslant\tilde{K}(u)$ and $\theta^p_{\{x,y\}}=1$. Proposition \ref{truncated negative level sets} is true for any choice of coupling probability $\tilde{\Q}^p$ satisfying \eqref{eq:Qdef}. Once its proof is completed, we will choose the coupling introduced in \eqref{isoeq}. This will automatically enforce the lower bound $-h$ required for the proof of Theorem \ref{percphitilde}, since $\tilde{\mathcal{I}}^u\subset\{x\in{\tilde{\Z}^d};\,\tilde{\phi}_x>-\sqrt{2u}\},$ see \eqref{htilde=0}. A good choice of $\theta^{p}_{\{x,y\}}$, cf.\ Lemma \ref{couplingwithbernoulli} below, will then allow to control the height of the field along the edges. To this effect, \eqref{A'u} essentially guarantees that on $\tilde{A}'_{u,p}$ we are dealing with Brownian bridges whose boundary values are uniformly bounded.

The proof of Proposition \ref{truncated negative level sets} follows a strategy very similar to the proof of Theorem~2.1 in \cite{MR3024098}, but we need to pay diligent attention to the dependence on $u$ in order to obtain the explicit bounds (\ref{Kpp'}). We use a renormalisation scheme akin to the one introduced in Section 4 of \cite{MR3024098}, which uses a sprinkling technique developed in \cite{MR2680403} and later improved in \cite{MR2891880} and \cite{MR3420516}. For $n\geqslant0$ and $L_0\geq 1$, we define the geometrically increasing sequence
\begin{equation}
\label{Ln}
    L_n=l_0^nL_0, \quad \text{where $l_0=4l(d)$ and $l(d)=4(5\cdot4^d+1)$}
\end{equation}
and the coarse-grained lattice model 
\begin{equation*}
    \mathbb{G}_0^{L_0}=L_0\Z^d \text{ and }\mathbb{G}_n^{L_0}=L_n\Z^d\subset\mathbb{G}_{n-1}^{L_0}\text{ for }n\geqslant1.
\end{equation*}
Note that, albeit only implicitly, the sequence $L_n$ depends on the choice of $L_0$, which is the only parameter in this scheme. For $x\in{\mathbb{G}_n^{L_0}},$ we further introduce the boxes
\begin{equation}
\label{Lambda}
    \Lambda_{x,n}^{L_0}=\mathbb{G}_{n-1}^{L_0}\cap(x+[0,L_n)^d),
\end{equation}
and note that $\{\Lambda_{x,n}^{L_0}; \, x\in{\mathbb{G}_n^{L_0}}\}$ forms a partition of $\mathbb{G}_{n-1}^{L_0}.$ For a given collection of events indexed by $\mathbb{G}_0^{L_0},$ that we denote by $A=(A_x)_{x\in{\mathbb{G}_0^{L_0}}},$ we define recursively the events $G_{x,n}^{L_0}(A)$ such that $G_{x,0}^{L_0}=A_x$ for all $x\in{\mathbb{G}_0^{L_0}},$ and for all $n\geqslant1$ and $x\in{\mathbb{G}_n^{L_0}},$
\begin{equation}
\label{GxnA}
    G_{x,n}^{L_0}(A) =\bigcup_{\stackrel{x_1,x_2\in{\Lambda_{x,n}^{L_0}}}{|x_1-x_2|_\infty\geqslant\frac{L_n}{l(d)}}}{G_{x_1,n-1}^{L_0}(A)\cap G_{x_2,n-1}^{L_0}(A)},
\end{equation}
where $|\cdot|_{\infty}$ stands for the $\ell^{\infty}$-distance on $\Z^d$. For each $x\in{\Z^d},$ let $T_x$ be the translation operator on the space of point measures on $W^*$, the space of doubly infinite trajectories on $\tilde{\Z}^d$ modulo time-shift such that, if $\mu$ is such a measure, then $T_x(\mu)$ is the point measure where each trajectory in the support of $\mu$ has been translated by $x$.  Moreover, in a slight abuse of notation, let $\tau_x$ be defined by
\begin{equation}
\label{eq:trans}
\tilde{\omega}\circ\tau_x=T_x(\tilde{\omega}).
\end{equation}
We introduce a family of events on the space $\Omega_{\text{coup}}$ on which $\tilde{\Q}^p$, cf.\ \eqref{eq:Qdef}, is defined. We say that an event $A\in{\sigma(\tilde{\phi}_z,\,z\in{\tilde{\Z}}^d})$ is increasing if there exists an increasing and measurable function $f:\Omega_0\rightarrow\{0,1\}$ such that, $\tilde{\Q}^p$-a.s., $1_A=f(\tilde{\phi}),$  and decreasing if $A^c$ is increasing. Recall that the events $E_{L_0}^{u,u'}$ and $F_{L_0}^{u,u'}$ from Definition \ref{defEandF} are respectively increasing and decreasing, and we will from now on tacitly consider them as subsets of $\Omega_{\text{coup}}.$

\begin{Def}
\label{good}
For each $u>0$, integer $L_0\geq 1$, $K>0$ and $p\in{[0,1]}$ let
\begin{enumerate}[(a)]
    \item $(\mathbf{E}_x^{L_0,u})_{x\in{\mathbb{G}_0^{L_0}}}$ be the family of increasing events such that, for all $x\in{\mathbb{G}_0^{L_0}},$ the event $\mathbf{E}_x^{L_0,u}=\tau_x^{-1}\left(E_{L_0}^{u,u}\right)$ occurs,
    \item $(\mathbf{F}_x^{L_0,u})_{x\in{\mathbb{G}_0^{L_0}}}$ be the family of decreasing events such that, for all $x\in{\mathbb{G}_0^{L_0}},$ the event $\mathbf{F}_x^{L_0,u}=\tau_x^{-1}\left(F_{L_0}^{u,u}\right)$ occurs,
\item $(\mathbf{C}_{x}^{L_0,K})_{x\in{\mathbb{G}_0^{L_0}}}$  be the family of decreasing events such that, for all $x\in{\mathbb{G}_0^{L_0}},$ the event $\mathbf{C}_{x}^{L_0,K}$ occurs if and only if for all $y\in{(x+[-1,2L_0+1)^d)\cap\Z^d},$ we have $\phi_y\leqslant K,$
\item $(\hat{\mathbf{C}}_x^{L_0,K})_{x\in{\mathbb{G}_0^{L_0}}}$  be the family of increasing events such that, for all $x\in{\mathbb{G}_0^{L_0}},$ the event $\hat{\mathbf{C}}_x^{L_0,K}$ occurs if and only if for all $y\in{(x+[-1,2L_0+1)^d)\cap\Z^d},$ we have $\phi_y\geqslant -K,$
\item $(\mathbf{D}_x^{L_0,p})_{x\in{\mathbb{G}_0^{L_0}}}$ be the family of events such that, for all $x\in{\mathbb{G}_0^{L_0}},$ the event $\mathbf{D}_{x}^{L_0,p}$ occurs if and only if for all $e\in{(x+[-1,2L_0+1)^d)\cap E},$ we have $\theta_e^{p}=1.$
\end{enumerate}
A vertex $x\in{\mathbb{G}_0^{L_0}}$ is called a good $(L_0,u,K,p)$ vertex if 
\begin{equation}
\label{goodevent}
    \mathbf{C}_x^{L_0,K}\cap\hat{\mathbf{C}}_x^{L_0,K}\cap \mathbf{D}_x^{L_0,p}\cap \mathbf{E}_x^{L_0,u}\cap \mathbf{F}_x^{L_0,u}
\end{equation}
occurs, and otherwise a bad $(L_0,u,K,p)$ vertex. 
\end{Def}
The reason for the choices in Definition \ref{defEandF} and \eqref{goodevent}, with regards to Proposition \ref{truncated negative level sets}, comes in the following.
\begin{Lemme}
\label{infinitepath} $(u>0, \, L_0 \geqslant 1, \, p\in(0,1])$

\medskip
\noindent If $(x_0,x_1,\dots)$ is an unbounded nearest neighbor path of good $(L_0,u,\tilde{K}(u),p)$ vertices in $\mathbb{G}_0^{L_0},$ then the set $\bigcup_{i=0}^\infty(x_i+[0,2L_0)^d)$ contains an unbounded connected path in $\tilde{A}'_{u,p},$ cf.~\eqref{A'u}.
\end{Lemme}
\begin{proof}
Let $x, y$ be two good $(L_0,u,\tilde{K}(u),p)$ vertices and neighbours in ${\mathbb{G}_0^{L_0}},$ and assume that there exists $e\in\{0,1\}^d$ such that $x+eL_0=y.$ Since $\mathbf{E}_{x}^{L_0,u}$ holds, there exist two random sets $\mathcal{A}_{x,0}\subset \tilde{\I}^u\cap(x+[0,L_0)^d)$ and $\mathcal{A}_{x,e}\subset \tilde{\I}^u\cap(x+eL_0+[0,L_0)^d)$ which are connected in $\tilde{\I}^u\cap(x+[0,2L_0)^d$), and such that the sum of the local times on the vertices of each of those two sets is larger than $\frac{3}{4}uL_0^d.$ Moreover, since $\mathbf{F}_{y}^{L_0,u}$ occurs, the sum of the local times on the vertices of $\mathcal{A}_{x,e}\cup\mathcal{A}_{y,0}$ is smaller than $\frac{5}{4}uL_0^d$ because $\mathcal{A}_{x,e}\cup\mathcal{A}_{y,0}\subset\tilde{\I}^u\cap(y+[0,L_0)^d).$ Hence, $\mathcal{A}_{x,e}\cap\mathcal{A}_{y,0}\neq\emptyset,$ and this implies that $\mathcal{A}_{x,0}$ is connected to $\mathcal{A}_{y,0}$ in $\tilde{\I}^u\cap(x+[0,2L_0)^d).$

Applying the above to each of the neighbors in our path $(x_0,x_1,\dots),$ we get that for all $i\in{\N_0},$ $\mathcal{A}_{x_i,0}$ is connected to $\mathcal{A}_{x_{i+1},0}$ in $\tilde{\I}^u\cap((x_i+[0,2L_0)^d)\cup(x_{i+1}+[0,2L_0)^d)).$ Thus, one can find an unbounded connected path in $\tilde{\I}^u\cap\bigcup_{i=0}^\infty(x_i+[0,2L_0)^d),$ and this path is actually in $\tilde{A}'_{u,p}$ since $|\phi_x|\leqslant \tilde{K}(u)$ and $\theta_e^{p}=1$ for all $x,e\in{\bigcup_{i=1}^\infty(x_i+[-1,2L_0+1)^d)}$ by Definition \ref{good}, (c), (d), (e). 
\end{proof}

To prove that an unbounded nearest neighbor path of good $(L_0,u,\tilde{K}(u),p)$ vertices in $\mathbb{G}_0^{L_0}$ exists for a suitable choice of the parameters, we pair our good (seed) events with the renormalisation scheme \eqref{GxnA} to show that, if being a good seed is typical, i.e., if it occurs with probability sufficiently close to $1,$ then the probability of being good ``at level $n$'' cf. \eqref{GxnA} and \eqref{nbad} below, is overwhelming. The respective bounds for all events of interest, cf.\ Definition \ref{good}, can be found in Lemmas \ref{G0NEandG0NF}, \ref{G0NC} and \ref{G0ND} below. We first consider the events $(\mathbf{E}_x^{L_0,u})_{x\in{\Z^d}}$ and $(\mathbf{F}_x^{L_0,u})_{x\in{\Z^d}},$ and take advantage of Corollary \ref{E0uandF0u} and Theorem \ref{couplingforinterlacements} to show the following. 
\begin{Lemme} $(u_0 > 0)$
\label{G0NEandG0NF}

\medskip
\noindent There exist $C_2=C_2(d,u_0)$ and $C'_2=C'_2(d,u_0)$ such that for all $u\in{(0,u_0]}$ and $L_0\geq 1$ with $L_0^{1/7}u\geqslant C_2,$
\begin{equation}
\label{G0NE}
    \tilde{\Q}^p\left[G_{0,n}^{L_0}\left((\mathbf{E}^{L_0,u})^c\right)\right]\leqslant2^{-2^n},
\end{equation}
and for all $u\in{(0,u_0]}$ and $L_0 \geq 1$ with $L_0^{1/7}u\geqslant C'_2,$
\begin{equation}
\label{G0NF}
    \tilde{\Q}^p\left[G_{0,n}^{L_0}\left((\mathbf{F}^{L_0,u})^c\right)\right]\leqslant2^{-2^n}.
\end{equation}
\end{Lemme}
\begin{proof}
We only prove (\ref{G0NE}). The proof of (\ref{G0NF}) is similar. Fix $\delta \in (0,1)$ as in Corollary~\ref{E0uandF0u}, and let $\delta'\in{(0,1)}$ be small enough such that
\begin{equation}
\label{u_n}
    u_n:=\frac{u(1-\delta)}{\prod_{k=1}^{n-1}\left(1-\frac{\delta'}{2^k}\right)}<u,\ \forall n\geqslant1
\end{equation}
(with $u_1= u(1-\delta)$). For all $x\in{\mathbb{G}_0^{L_0}},$ let $\mathbf{E}_x^{L_0,u,u'}=\tau_x^{-1}(E_{L_0}^{u,u'}),$ cf. \eqref{eq:trans} and Definition~\ref{defEandF}, and note that $\mathbf{E}_x^{L_0,u,u}=\mathbf{E}_x^{L_0,u}.$ For all $u\in{(0,u_0]},$ $L_0\geqslant1,$ positive integer $n,$ $i\in{\{1,2\}}$ and $x_i\in{\mathbb{G}_n^{L_0}}$ such that $|x_1-x_2|_\infty\geqslant\frac{L_{n+1}}{l(d)}=4L_n,$ the events $G_{x_i,n-1}^{L_0}\left((\mathbf{E}^{L_0,u',u})^c\right)$ are $\sigma(\tilde{\ell}_{z,u},\  z\in{x_i+[0,L_n+L_0)^d})$ measurable for all $u'>0$, and, defining $r_n:=\frac{L_n+L_0}{2},$
\begin{equation*}
s_n:=d\left(x_1+[0,2r_n)^d,x_2+[0,2r_n)^d]\right)\geqslant d(x_1,x_2)-r_n\geqslant L_n.   
\end{equation*}
By Theorem \ref{couplingforinterlacements} applied with $\eps=1-\frac{u_n}{u_{n+1}}=\delta'2^{-n},$ and since the events $(\mathbf{E}^{L_0,u',u})^c$ are decreasing, there exist two constants $C$ and $c$ independent of $u,$ $n$ and $L_0$ such that 
\begin{equation}
\label{decoupapp}
\begin{split}
    &\tilde{\Q}^p\left[G_{x_1,n-1}^{L_0}\left((\mathbf{E}^{L_0,u_{n+1},u})^c\right)\cap G_{x_2,n-1}^{L_0}\left((\mathbf{E}^{L_0,u_{n+1},u})^c\right)\right] \\
    &\  \leqslant\tilde{\Q}^p\left[G_{x_1,n-1}^{L_0}\left((\mathbf{E}^{L_0,u_n,u})^c\right)\right]\tilde{\Q}^p\left[G_{x_2,n-1}^{L_0}\left((\mathbf{E}^{L_0,u_n,u})^c\right)\right]
    +C(s_n+r_n)^d\exp\left(-cus_n^{d-2}4^{-n}\right).
\end{split}
\end{equation}
We have chosen $l_0^{d-2}>8,$ see \eqref{Ln}, whence for $L_0^{1/7}u\geqslant c(d,u_0),$
\begin{equation}
\label{decoupapp2}
    l_0^{2d}C(s_n+L_n+L_0)^d\exp\left(-cus_n^{d-2}4^{-n}\right)\leqslant C'\exp\big(-cuL_0^{1/7}2^n\big)\leqslant\frac{1}{(4l_0^{2d})2^{2^{n+1}}}.
\end{equation}
We now prove by induction over $n$ that for all $x\in{\mathbb{G}_n^{L_0}},$ and all $u\in{(0,u_0]},$
\begin{equation}
\label{induction}
    \tilde{\Q}^p\left[G_{x,n}^{L_0}\left((\mathbf{E}^{L_0,u_{n+1},u})^c\right)\right]\leqslant\frac{1}{(2l_0^{2d})2^{2^n}},  \quad \text{ if } L_0^{1/7}u \geq c(d,u_0).
\end{equation}
For $n=0,$ the bound on the right-hand side of \eqref{induction} is purely numerical. Thus, it is clear from Corollary \ref{E0uandF0u}, and since $G_{x,0}^{L_0}\left((\mathbf{E}^{L_0,u_1,u})^c\right)=(\mathbf{E}_x^{L_0,u_1,u})^c,$ see above \eqref{GxnA}, that if one takes $L_0^{1/7}u$ large enough (only depending on $u_0$ and $d$), then (\ref{induction}) holds for $n=0$ on account of \eqref{E0u}. Suppose now it holds for $n-1 \geq 0$. Then, according to (\ref{GxnA})
\begin{align*}
    \tilde{\Q}^p\left[G_{x,n}^{L_0}\left((\mathbf{E}^{L_0,u_{n+1},u})^c\right)\right]
&\leqslant \sum_{\stackrel{x_1,x_2\in{\Lambda_{x,n}^{L_0}}}{|x_1-x_2|\geqslant \frac{L_{n}}{l(d)}}}\tilde{\Q}^p\left[G_{x_1,n-1}^{L_0}\left((\mathbf{E}^{L_0,u_{n+1},u})^c\right)\cap G_{x_2,n-1}^{L_0}\left((\mathbf{E}^{L_0,u_{n+1},u})^c\right)\right]
    \\&\leqslant\frac{1}{(2l_0^{2d})2^{2^{n}}},
\end{align*}
where the last equality follows from \eqref{decoupapp}, \eqref{decoupapp2}, the induction hypothesis and $|\Lambda_{x,n}^{L_0}|\leqslant l_0^d,$ and (\ref{induction}) follows. The claim (\ref{G0NE}) then follows from \eqref{u_n}, \eqref{induction} and the fact that the $\left(\mathbf{E}^{L_0,u',u}\right)^c$ are decreasing events.
\end{proof}

We now turn to the Gaussian free field part $\mathbf{C}^{L_0,K}$ and $\hat{\mathbf{C}}^{L_0,K},$ see (c) and (d) in Definition \ref{good}, of the good events in \eqref{goodevent}. Sprinkling techniques have been used successfully in investigating  level set percolation of the Gaussian free field, see for example \cite{MR3053773}, \cite{MR3163210} or \cite{MR3339867}, and also \cite{2016arXiv161202385R} with regard to non-Gaussian measures. These techniques imply similar results as for random interlacements, and this is mainly due to the fact that decoupling inequalities as \eqref{decinc} also hold for the Gaussian free field. With hopefully obvious notation, in writing $\Phi + c$ below, with $\Phi$ as in \eqref{Phi}, we mean the field whose value is shifted by $c \in \mathbb{R}$ everywhere. 

\begin{The}[{\cite[Corollary 1.3]{MR3325312}} and thereafter]
\label{couplingforGFF}
Let $A_1$ and $A_2$ be two non intersecting subsets of $\Z^d$, define $s=d(A_1,A_2)$ and assume that the minimum $r$ of their diameters is finite. Then, there exist positive constants $\kappa'_0(d)$ and $\kappa'_1(d)$ such that, for all $\epsilon\in{(0,1)}$, and any two functions $f_i:\R^{\Z^d}\rightarrow[0,1]$ which are $\sigma(\Phi_x,\, x\in{A_i})$ measurable for each $i\in{\{1,2\}},$ and either both increasing or both decreasing,
   \begin{equation}
    \label{increasinggff}
    \mathbb{E}^G[f_1(\Phi)f_2(\Phi)]\leqslant \mathbb{E}^G[f_1(\Phi\pm\eps)]\mathbb{E}^G[f_2(\Phi\pm\eps)]+\kappa'_0(r+s)^d\exp(-\kappa'_1\epsilon^2s^{d-2}),
    \end{equation}
where the plus sign corresponds to the case where the $f_i$'s are increasing and the minus sign to the case where the $f_i$'s are decreasing.
 
\end{The}
Theorem 1.2 in \cite{MR3325312} gives a slightly better inequality, but \eqref{increasinggff} will be sufficient for our purposes, and readily yields the following analogue of Lemma \ref{G0NEandG0NF} for the events pertaining to the free field. 
\begin{Lemme}
\label{G0NC}
There exist constants $C_3(d)>1$ and $C'_3(d)>0$ such that for all $L_0\geqslant C_3$ and $K>0$ with 
\begin{equation}
\label{C'3}
K\geqslant C'_3\sqrt{\log(L_0)},
\end{equation}
one has
\begin{equation}
\label{G0NCeq}
    \tilde{\Q}^p\left[G_{0,n}^{L_0}\Big((\mathbf{C}^{L_0,K})^c\Big)\right]\leqslant2^{-2^n}\quad \text{ and } \quad \tilde{\Q}^p\left[G_{0,n}^{L_0}\left((\hat{\mathbf{C}}^{L_0,K})^c\right)\right]\leqslant2^{-2^n}.
\end{equation}
\end{Lemme}
\begin{proof}
One knows from (2.35) and (2.38) in \cite{MR3053773} that if $K>C\sqrt{\log(L_0)}$ for some constant $C$ large enough,
\begin{equation}
\label{C0K}
     \tilde{\Q}^p\left[\left(\hat{\mathbf{C}}_0^{L_0,K}\right)^c\right]=\tilde{\Q}^p\left[\Big(\mathbf{C}_0^{L_0,K}\Big)^c\right]=\tilde{\Q}^p\Big(\sup_{x\in{[-1,2L_0+1)^d}}{\phi_x}>K\Big)\leqslant e^{-\frac{\left(K-C\sqrt{\log(L_0)}\right)^2}{2g(0)}}.
\end{equation}
The claim (\ref{G0NCeq}) now follows by induction over $n$ from (\ref{C0K}) and Theorem \ref{couplingforGFF} in exactly the same way as Lemma \ref{G0NEandG0NF} was obtained from Corollary \ref{E0uandF0u} and Theorem \ref{couplingforinterlacements}.
\end{proof}

Finally, we collect a simple estimate for the Bernoulli part of our good events $\mathbf{D}^{L_0,p},$ see part (e) in Definition \ref{good}.
\begin{Lemme}[{\cite[Lemma 4.7]{MR3024098}}]
\label{G0ND}
There exists $C_4=C_4(d)$ such that for all $L_0\geqslant1$ and $p\in{(0,1)}$ satisfying
\begin{equation}
\label{inequalityp}
    p\geqslant\exp\left(-\frac{C_4}{L_0^d}\right),
\end{equation}
one has\begin{equation}
\label{inequalityp2}
    \tilde{\Q}^p\left[G_{0,n}^{L_0}\left((\mathbf{D}^{L_0,p})^c\right)\right]\leqslant2^{-2^n}.
\end{equation}
\end{Lemme}

The bounds of Lemmas \ref{G0NEandG0NF}, \ref{G0NC} and \ref{G0ND} allow for a proof of Proposition \ref{truncated negative level sets} by means of a standard duality argument. On account of Lemma \ref{infinitepath}, this requires an estimate on the probability to see certain long (dual) paths. The relevant events, see \eqref{eq:Hevent}, can be suitably expressed in terms of bad vertices at level $n$, as Lemma \ref{Hinclbad} asserts. 

Recall the definition of good $(L_0,u,K,p)$ vertices in \eqref{goodevent}. For $n\geq 0$, we call $x\in{\mathbb{G}_n^{L_0}}$ a bad $n-(L_0,u,K,p)$ vertex if the event
\begin{equation}
\label{nbad}
    G^{L_0}_{x,n}\big((\mathbf{C}^{L_0,K})^c\big)\cup G^{L_0}_{x,n}\big((\hat{\mathbf{C}}^{L_0,K})^c\big)\cup G^{L_0}_{x,n}\big( ( \mathbf{D}^{L_0,p})^c\big)\cup G^{L_0}_{x,n}\big((\mathbf{E}^{L_0,u})^c\big)\cup G^{L_0}_{x,n}\big((\mathbf{F}^{L_0,u})^c\big)
\end{equation}
occurs, and a good $n-(L_0,u,K,p)$ vertex otherwise. Note that a good $0-(L_0,u,K,p)$ vertex is simply a good $(L_0,u,K,p)$ vertex. We say that $(x_0,x_1,\dots,x_n,\dots)$ is a $*$-path in $\mathbb{G}_0^{L_0}$ if for all $i\in \{ 0,1,\dots \}$,  $x_i\in{\mathbb{G}_0^{L_0}}$ and $\|x_i-x_{i+1}\|_{\infty}=L_0.$ For each $u,K>0,$ integer $L_0\geq 1,$ $p\in{(0,1)},$ $0<M<N$ with $M,\ N$ multiples of $L_0,$ and $x\in{\mathbb{G}_0^{L_0}},$ let 
\begin{equation}
\begin{split}
\label{eq:Hevent}
H_M^N(x;L_0,u,K,p) = \big\{&\text{$(x+[-M,M]^d)$ is connected to $(x+\partial[-N,N]^d)$ by}\\
&\text{a $*$-path of bad $(L_0,u,K,p)$ vertices in $\mathbb{G}_0^{L_0}$}\big\}.
\end{split}
\end{equation}
 Here, $\partial[-N,N]^d$ denotes the boundary of the set $[-N,N]^d,$ which intersects $\mathbb{G}_0^{L_0}$ since $N$ is a multiple of $L_0.$ The following lemma asserts that $H_{L_n}^{2L_n}(x;L_0,u,K,p)$ can only happen if there is a bad $n-(L_0,u,K,p)$ vertex in the box of radius $2L_n$ around $x.$
\begin{Lemme}
\label{Hinclbad}
For all integers $n\geq 0$ and $L_0\geq 1,$ $u,K>0,$  $p\in{(0,1)},$ and $x\in{\mathbb{G}_n^{L_0}},$ 
\begin{equation}
\label{incl}
    H_{L_n}^{2L_n}(x;L_0,u,K,p)\subset\bigcup_{y\in{\mathbb{G}_n^{L_0}\cap(x+[-2L_n,2L_n)^d)}}\{y\text{ is }n-(L_0,u,K,p) \text{ bad}\}.
\end{equation}
\end{Lemme}
\begin{proof}
This is a consequence of Lemma 4.4 of \cite{MR3163210} (with $N\equiv5$, $r\equiv l(d)$ and $L_0, l_0$ as in \eqref{Ln} above). We include the proof for the reader's convenience. We proceed by  induction over $n$: it is clear that \eqref{incl} is true for $n=0,$ and we assume that it holds for any choice of $x$ up to level $n-1.$ If $H_{L_n}^{2L_n}(x;L_0,u,K,p)$ occurs, there exists a $*$-path $\pi$ of bad $(L_0,u,K,p)$ vertices in $\mathbb{G}_0^{L_0}$ from $(x+[-L_n,L_n]^d)$ to $(x+\partial[-2L_{n},2L_n]^d).$ This path intersects the concentric $\ell^{\infty}$-spheres  $(x+\partial[-L_n-16iL_{n-1},L_n+16iL_{n-1}]^d)$ for all $i\in{\left\{0,\dots,m-1\right\}},$ where $m=5\cdot4^d+1$ (recall that $l_0=16m$). In view of \eqref{Ln}, for all $i\in{\left\{0,\dots,m-1\right\}},$ one can thus find $y_i\in{\mathbb{G}_{n-1}^{L_0}}\cap(x+\partial[-L_n-16iL_{n-1},L_n+16iL_{n-1}]^d)$ such that $\pi\cap(y_i+[-L_{n-1},L_{n-1}]^d)\neq\emptyset.$ 

For each $i\in{\left\{0,\dots,m-1\right\}},$ since $(y_i+[-2L_{n-1},2L_{n-1}]^d)\subset(x+[-2L_n,2L_n]^d),$ the connected $*$-path $\pi$ in $\mathbb{G}_0^{L_0}$  connects $(y_i+[-L_{n-1},L_{n-1}]^d)$ to $(y_i+\partial[-2L_{n-1},2L_{n-1}]^d),$ and thus the induction hypothesis implies that there exists $z_i\in{(y_i+[-2L_{n-1},2L_{n-1})^d)}$ which is $(n-1)-(L_0,u,K,p)$ bad, and in particular $z_i\in{\mathbb{G}_{n-1}^{L_0}}.$ There are $m=5\cdot4^d+1$ such $z_i,$ and since there are only $4^d$ elements in $\mathbb{G}_n^{L_0}\cap\left(x+[-2L_n,2L_n)^d\right),$ one can find $x_0$ in this set such that $\Lambda_{x_0,n}^{L_0},$ cf.\ \eqref{Lambda}, contains at least $6$ different $z_i.$ By \eqref{nbad}, one can thus find $k\neq j$ in  $\left\{0,\dots,m-1\right\}$ and $A_0\in{\{(\mathbf{C}^{L_0,K})^c,(\hat{\mathbf{C}}^{L_0,K})^c,(\mathbf{D}^{L_0,p})^c,(\mathbf{E}^{L_0,u})^c,(\mathbf{F}^{L_0,u})^c\}}$ such that $z_k$ and $z_j$ are in $\Lambda_{x_0,n}^{L_0},$ and $G_{z_k,n-1}^{L_0}(A_0)$ and $G_{z_j,n-1}^{L_0}(A_0)$ both occur. Moreover,
\begin{equation*}
    \|z_k-z_j\|_{\infty}\geqslant\|y_k-y_j\|_{\infty}-4L_{n-1}\geqslant12L_{n-1}\geqslant {L_n}/{l(d)},
\end{equation*}
which, in view of \eqref{GxnA}, implies that $G_{x_0,n}^{L_0}(A_0)$ occurs, and thus $x_0$ is $n-(L_0,u,K,p)$ bad.
\end{proof}

By Lemmas \ref{G0NEandG0NF}, \ref{G0NC} and \ref{G0ND}, we know that for all $u\in{(0,u_0]},$ and for a suitable choice of the parameters $L_0,\ K$ and $p,$ the probability that a vertex is $n-(L_0,u,K,p)$ bad is very small. Lemma \ref{Hinclbad} then yields that a $*$-path of $(L_0,u,K,p)$ bad vertices in $\mathbb{G}_0^{L_0}$ exists with very small probability only, and on account of Lemma \ref{infinitepath}, we can prove Proposition~\ref{truncated negative level sets} using a Peierls argument. 

\begin{proof}[Proof of Proposition \ref{truncated negative level sets}]
Choose a constant $C_5=C_5(d,u_0)$ large enough such that, upon defining
\begin{equation*}
    L_0(u)=\left\lceil C_5 / u^7\right\rceil,
\end{equation*}
one has $L_0(u)^{1/7}u\geqslant \max(C_2,C'_2),$ cf.\ Lemma \ref{G0NEandG0NF}, and $L_0(u)\geqslant C_3,$ cf.\ Lemma~\ref{G0NC}, for all $u\in{(0,u_0]}.$ One can now find constants $c_1,\ C_1,\ c'_1\text{ and } C'_1$ such that if (\ref{Kpp'}) holds, then $\tilde{K}(u)\geqslant C'_3\sqrt{\log(L_0(u))},$ cf.\ \eqref{C'3}, and $p(u)\geqslant\exp(-\frac{C_4}{L_0(u)^d}),$ cf.\ \eqref{inequalityp}, for all $u\in{(0,u_0]}.$ Let us now fix arbitrarily some $u\in{(0,u_0]}$ and $p\in[p(u),1].$ Lemma \ref{Hinclbad}, Lemmas~\ref{G0NEandG0NF} and \ref{G0NC} and Lemma \ref{G0ND} can now be applied with $L_0=L_0(u),$ $K=\tilde{K}(u)$ and $p$, to yield
\begin{IEEEeqnarray*}{rCl}
\IEEEeqnarraymulticol{3}{l}{\tilde{\Q}^{p}\left(H_{L_n}^{2L_n}(0;L_0(u),u,\tilde{K}(u),p)\right)}
\\\ \ &\stackrel{\eqref{incl}}{\leqslant}&4^d\tilde{\Q}^{p}\big(0\text{ is }n-(L_0(u),u,\tilde{K}(u),p)\text{ bad}\big)
\\&\stackrel{\eqref{nbad}}{\leqslant}&4^d\Big\{\tilde{\Q}^{p}\left[G_{0,n}^{L_0(u)}\left(\big(\mathbf{C}^{L_0(u),\tilde{K}(u)}\big)^c\right)\right]+\tilde{\Q}^{p}\left[G_{0,n}^{L_0(u)}\left(\big(\hat{\mathbf{C}}^{L_0(u),\tilde{K}(u)}\big)^c\right)\right]
\\&\phantom{aa}+&\tilde{\Q}^{p}\left[G_{0,n}^{L_0(u)}\left(\left(\mathbf{D}^{L_0(u),p}\right)^c\right)\right]+\tilde{\Q}^{p}\left[G_{0,n}^{L_0(u)}\left(\left(\mathbf{E}^{L_0(u),u}\right)^c\right)\right]+\tilde{\Q}^{p}\left[G_{0,n}^{L_0(u)}\left(\left(\mathbf{F}^{L_0(u),u}\right)^c\right)\right]\Big\}
\\&\leqslant&5\cdot4^d\cdot2^{-2^n},
\end{IEEEeqnarray*}
using \eqref{G0NE}, \eqref{G0NF}, \eqref{G0NCeq} and \eqref{inequalityp2} in the last step.
Since this bound holds for all $n\geqslant 0,$ and, in view of \eqref{Ln}, $H_0^N(0;L_0(u),u,\tilde{K}(u),p)\subset H_{L_n}^{2L_n}(0;L_0(u),u,\tilde{K}(u),p)$ for any $n\in{\N}$ such that $2L_n\leqslant N,$ one can find constants $c,C>0$ depending only on $d,$ $u$ and $u_0$ such that, for all integers $N,$
\begin{equation}
\label{H0N}
    \tilde{\Q}^{p}\left(H_{0}^{N}\big(0;L_0(u),u,\tilde{K}(u),p\big)\right)\leqslant C\exp\left(-cN^c\right).
\end{equation}
Given \eqref{H0N}, the argument proceeds as follows. For any set $A\subset\Z^2\times\{0\}^{d-2},$ define $(x_0,\dots,x_n,\dots)$ to be a nearest neighbor path of good $(L_0(u),u,\tilde{K}(u),p)$ vertices in $\mathbb{G}_0^{L_0(u)}\cap(\Z^2\times\{0\}^{d-2})$ that connects $A$ to $\infty$ if the $x_i\in{\mathbb{G}_0^{L_0(u)}\cap(\Z^2\times\{0\}^{d-2})}$ are good $(L_0(u),u,\tilde{K}(u),p)$ vertices, $\|x_i-x_{i+1}\|_1=L_0(u)$ for all $i\in\{ 0,1, \dots\}$, $x_0\in{A}$ and $\|x_i\|_{\infty}\rightarrow\infty,$ as $i\rightarrow\infty.$ Now, assume that there exists no unbounded nearest neighbor path of good $(L_0(u),u,\tilde{K}(u),p)$ vertices in $\mathbb{G}_0^{L_0(u)}\cap(\Z^2\times\{0\}^{d-2}),$ and in particular that for all $M\in{L_0(u) \cdot\N} \equiv \mathbb{N}_u$ there is no nearest neighbor path of good $(L_0(u),u,\tilde{K}(u),p)$ vertices that connects $[M,M]^2\times\{0\}^{d-2}$ to $\infty.$ Then by planar duality, for all $M\in \mathbb{N}_u,$ there exists a $*$-path $\pi$ around $[-M,M]^2\times\{0\}^{d-2}$ in $\mathbb{G}_0^{L_0(u)}\cap(\Z^2\times\{0\}^{d-2})$ of bad $(L_0(u),u,\tilde{K}(u),p)$ vertices. If $N\geqslant M$ denotes the smallest multiple of $L_0(u)$ such that $x_N\equiv(N,0)\in{\N_u\times\{0\}^{d-1}}$ is in $\pi,$ then $H_0^N(x_N;L_0(u),u,\tilde{K}(u),p)$ occurs. Thus, the probability that there is no infinite nearest neighbor path of good $(L_0(u),u,\tilde{K}(u),p)$ vertices in $\mathbb{G}_0^{L_0(u)}\cap(\Z^2\times\{0\}^{d-2})$ that connect $[-M,M]^2\times\{0\}^{d-2}$ to $\infty$ is bounded by
\begin{equation*}
    \sum_{N\in \mathbb{N}_u: \, N\geqslant M}{\tilde{\Q}^{p}\left(H_{0}^{N}(0;L_0(u),u,\tilde{K}(u),p\right)}\leqslant    \sum_{N\in \mathbb{N}_u: \, N\geqslant M}{C\exp\left(-cN^c\right)}.
\end{equation*}
This is true for all $M\in \mathbb{N}_u,$ hence the probability of having no unbounded nearest neighbor path of good $(L_0(u),u,\tilde{K}(u),p)$ vertices in $\mathbb{G}_0^{L_0(u)}\cap(\Z^2\times\{0\}^{d-2})$ is $0.$ Lemma \ref{infinitepath} then implies that the set $\tilde{A}'_{u,p}$ percolates (almost surely), for any $u\in{(0,u_0]}$ and $p\in{[p(u),1]}$, and
the claim of Proposition \ref{truncated negative level sets} follows.
\end{proof}
With Proposition \ref{truncated negative level sets} at hand, it is possible to deduce Theorem \ref{percphitilde} for a good choice of coupling $\tilde{\Q}^p$ in \eqref{eq:Qdef}. The idea is to use \eqref{isoeq}, and to suitably couple the Bernoulli percolation $\mathcal{B}^{p}$ with $\{|\tilde{\phi}|\leqslant K(h)\}$ on the edges. 

The key to the proof of Theorem \ref{percphitilde} is the following lemma, by which one can essentially couple a Bernoulli percolation $\{\mathcal{B}^{p}=1\}$ on the edges with sufficiently large success parameter $p\geqslant 1-C'_1u^{c'_1},$ cf.\ \eqref{Kpp'}, with $\{|\tilde{\phi}|\leqslant K'\}$ on the edges for $K'$ large enough. 

\begin{Lemme}
\label{couplingwithbernoulli}
Let $\tilde{\phi}$ be a Gaussian free field on the cable system under $\tilde{\P}.$ For all $u_0>0,$ there exist positive constants $C_0$ and $c_0$ such that, for all $u\in{(0,u_0]},$ 
with $\tilde{K}(u)$ and $p(u)$ as defined in \eqref{Kpp'}, $h=\sqrt{2u}$ and $K(h)$ as defined in \eqref{K}, the following holds: under $\tilde{\P},$ there exists a family of independent Bernoulli variables $\mathcal{B}^{\tilde{p}(u)}=(\theta_e^{\tilde{p}(u)})_{e\in{E}}$ with parameter $\tilde{p}(u) \geqslant p(u),$ and the property that
\begin{equation}
\begin{split}
    \label{couplingontheedges}
    &\text{for all } e=\{x,y\}\in{E}\text{, if }|\phi_x|\leqslant\tilde{K}(u)\text{ and }|\phi_y|\leqslant\tilde{K}(u),\\
    &\text{then } \big\{ \theta_{e}^{{\tilde{p}(u)}}=1\Rightarrow\forall\,z\in{\overline{I_e}},\ |\tilde{\phi}_z|\leqslant K(h) \big\}.
    \end{split}
\end{equation}
\end{Lemme}

\begin{proof}
Let $u\in{(0,u_0]}$ and $h=\sqrt{2u}.$  With $C_1,$ $c_1,$ $C'_1$ and $c'_1$ as given by Proposition \ref{truncated negative level sets}, fix constants $C_0$ and $c_0$ depending only on $u_0$ and $d$ such that 
\begin{equation}
\begin{split}
\label{eq:C_0fix}
  K(h)\stackrel{\eqref{K}}{ =}  \sqrt{\log\left(\frac{C_0}{(2u)^{c_0/2}}\right)}
  &\geqslant\sqrt{\log\left(\frac{C_1}{u^{c_1}}\right)}+\sqrt{-\frac{1}{2}\log\left(\frac{C'_1u^{c'_1}}{2}\right)} \\
  &\stackrel{ \eqref{Kpp'}}{=} \tilde{K}(u) + \sqrt{-\frac{1}{2}\log\left(\frac{1-p(u)}{2}\right)}.
\end{split}
\end{equation}
Let $(B^e)_{e\in{E}}$ be defined as in \eqref{brownianbridgeontheedges}, and recall that $(B^e)_{e\in{E}}$ is an i.i.d.\ family of Brownian bridges with length $\frac{1}{2}$ of a Brownian motion with variance $2$ at time $1.$ For all $e\in{E}$, define
\begin{equation}
\label{definitionofthetatilde}
    \theta_e^{{\tilde{p}(u)}}=\left\{
    \begin{array}{ll}
        1,&\text{ if } |B_t^e|\leqslant K(h)-\tilde{K}(u) \text{ for all }  t\in{\left[0,\frac{1}{2}\right],}\\
        0,& \text{ otherwise.}
    \end{array}
    \right.
\end{equation}
Then $\big(\theta_e^{\tilde{p}(u)}\big)_{e\in{E}}$ is an i.i.d.\ family of Bernoulli variables with parameter
\begin{equation*}
    {\tilde{p}(u)}\stackrel{\text{def.}}{=}\tilde{\P}\left(\forall\,t\in{[0,1/2]},\ |B_t^e|\leqslant K(h)-\tilde{K}(u)\right).
\end{equation*}
Moreover, by symmetry (the boundary values of $B^e$ are both $0$) and Lemma \ref{Brownianbridge},
\begin{equation*}
    {\tilde{p}(u)}\geqslant1-2\tilde{\P}\left(\sup_{t\in{\left[0,1/2\right]}}B_t^e\geqslant K(h)-\tilde{K}(u)\right)
    \stackrel{\eqref{sup}}{\geqslant}1-2\exp\left(-2(K(h)-\tilde{K}(u))^2\right)
     \stackrel{\eqref{eq:C_0fix}}{\geqslant}p(u),
\end{equation*}
and, using (\ref{brownianbridgeontheedges}) and (\ref{definitionofthetatilde}), for all $e=\{x,y \}\in{E}$ such that $\phi_{x}\leqslant \tilde{K}(u)$ and $\phi_{y}\leqslant \tilde{K}(u),$
\begin{equation*}
\begin{split}
    \theta_e^{{\tilde{p}(u)}}=1&\quad \Rightarrow \quad \forall t\in{\left[0,1/2\right]},\ |\tilde{\phi}_{x+tv_{(x,y)}}-(1-2t)\phi_{x}-2t\phi_{y}|\leqslant K(h)-\tilde{K}(u)\\
    & \quad \Rightarrow \quad \forall z\in{\overline{I_e}},\ 
    |\tilde{\phi}_{z}|\leqslant K(h),
    \end{split}
\end{equation*}
whence \eqref{couplingontheedges}.
\end{proof}

\begin{proof}[Proof of Theorem \ref{percphitilde}]
Let $u_0=\frac{h_0^2}{2},$ and, for any $h\in{(0,h_0]}$, define $u=\frac{h^2}{2}.$ Let $\tilde{\P}^u$ be the coupling from Theorem \ref{iso}, under which there exist a Gaussian free field $\tilde{\phi}$ and a random interlacement process $\tilde{\omega}$ such that \eqref{htilde=0} holds. For this $\tilde{\phi}$, let $\mathcal{B}^{\tilde{p}(u)} = \mathcal{B}^{\tilde{p}(u)}(\tilde{\phi}) $ be the family of independent Bernoulli variables under $\tilde{\P}^u$ introduced in Lemma \ref{couplingwithbernoulli}. This yields a coupling $\tilde{\Q}^{\tilde{p}(u)}$ satisfying \eqref{eq:Qdef}, with parameter $ \tilde{p}(u)\geqslant p(u).$ One can now apply Proposition \ref{truncated negative level sets} to obtain that, $\tilde{\P}^u$-a.s, the set $\tilde{A}'_{u,\tilde{p}(u)},$ cf.\ \eqref{A'u}, contains an unbounded connected component in the thick slab $\tilde{\Z}^2\times[0,L_0(u))^{d-2}$, and thus \eqref{couplingontheedges} yields that $\tilde{\P}^u$-a.s.\ the set
\begin{equation}
\label{percolationforIu}
    \big(\tilde{\I}^u\cap(\tilde{\Z}^d\setminus{\Z^d})\big)\cup\left\{x\in{\Z^d}; \ \forall\,v\in{V^0},\ \forall t\in{\left[0,1/2\right]},\ |\tilde{\phi}_{x+tv}|\leqslant K(h)\right\}
\end{equation}
contains an unbounded connected component in the thick slab $\tilde{\Z}^2\times[0,L_0(u))^{d-2}.$ Now \eqref{htilde=0} implies that the set defined in (\ref{percolationforIu}) is included in $\tilde{A}_h(\tilde{\phi})$, and Theorem \ref{percphitilde} follows.
\end{proof}

\section{Percolation for positive level set}
\label{proofofthemainresult}

In this section, we prove our main result, Theorem \ref{mainresult}, with the help of Theorem \ref{percphitilde}. We consider the Gaussian free field $\tilde{\Phi}$ on $\tilde{\Z}^d$ as defined in  \eqref{tildephi}, and, with Theorem \ref{percphitilde} at hand, we will no longer need random interlacements nor the coupling \eqref{isoeq} to prove Theorem \ref{mainresult}. A key ingredient is the following observation: we have shown, see \eqref{1:30} that the set $\{x\in{\tilde{\Z}^d};\,-h\leqslant\tilde{\Phi}_x\leqslant K(h)\}$ contains an unbounded connected component for large enough $K(h)$, cf.\ \eqref{K}. 
Suppose that $x \in \Z^d$ is a vertex inside this unbounded component, and that $I_e$ is attached to $x$ (recall that $\Phi =\tilde \Phi \restriction \Z^d$). Then, since $\tilde{\Phi}$ behaves like a Brownian bridge on $\overline{I_e},$ see \eqref{brownianbridgeontheedges}, the probability that $\tilde{\Phi}_z\geqslant-h$ for all $z\in{I_e}$ becomes very small as $h \searrow 0$. In fact, since $hK(h)\rightarrow0$ as $h\searrow0,$ if $e=\{x,y\},$ for sufficiently small $h>0$, it is more costly to keep $\tilde{\Phi}\geqslant-h$ along the entire cable $\overline{I_e},$ than to require $\Phi_x\geqslant h$ (\textit{at} the vertex $x$ only!), knowing that $-h\leqslant\Phi_x\leqslant K(h)$ and $|\Phi_y|\leqslant K(h)$, see Lemma \ref{comparingEandF} for the corresponding statement. Accordingly, the probability that the set $\{x\in{\tilde{\Z}^d};\, -h\leqslant\tilde{\Phi}_x\leqslant K(h)\}$ contains an unbounded connected component becomes smaller than the probability that the set $\{x\in{\Z^d};\, h\leqslant\Phi_x\leqslant K(h)\}$ contains an infinite cluster (in $\Z^d$) as $h$ goes to $0,$ which implies Theorem \ref{mainresult}.

Comparing the probability that $\Phi_x\geqslant h$ knowing that $-h\leqslant\Phi_x\leqslant K(h)$ with the probability that the Brownian bridge on $I_e$ remains above level $-h$ in a uniform way requires some control on the Gaussian free field $\tilde\Phi$ in the neighborhood of $x$, and for this purpose we are actually going to use Theorem \ref{percphitilde} and not only \eqref{1:30}. We define, for $x\in{\Z^d}$ and $v\in{V^0},$ the subsets $U^{x,v}$ and $U^x$ of $\tilde{\Z}^d$ by 
\begin{equation}
\label{eq:Usets}
U^{x,v}=x+\Big[0,\frac{1}{4}v\Big)\text{ and }U^x=\bigcup_{v\in{V^0}}U^{x,v}=x+\bigcup_{v\in{V^0}}\Big[0,\frac{1}{4}v\Big).
\end{equation}
We call $\mathcal{K}^x\equiv\partial U^x$ the boundary of $U^x,$ which has exactly $2d$ elements, and define $\mathcal{K}=\bigcup_{x\in{\Z^d}}\mathcal{K}^x.$ Henceforth, we set
\begin{equation}
\label{h_0choice}
    h_0=1
\end{equation}
in all the previous definitions and results, and in particular in Theorem \ref{percphitilde} (this value is chosen arbitrarily in $(0,\infty)$). For any $h\in{(0,1]},$ we define $K(h)$ as in \eqref{K} (with $c_0, \, C_0$ numerical constants depending only on $d$ by the choice \eqref{h_0choice}). We further define two families of events $(E_{h}^{x,v})_{x\in{\Z^d}}$ and $(F_{h}^{x,v})_{x\in{\Z^d},v\in{V^0}}$ (part of $\Omega_0$, cf. above \eqref{tildephi}) by
\begin{equation}
\label{FFFF}
E_{h}^{x,v}=\left\{\tilde{\Phi}_{x+\frac{1}{4}v}\geqslant-h\right\}\cap\Big\{\forall y\in{\mathcal{K}^x};\, |\tilde{\Phi}_y|\leqslant K(h)\Big\}\text{ and }
F_{h}^{x,v}=\left\{\forall\,z\in{U^{x,v}};\, \tilde{\Phi}_{z}\geqslant -h\right\},
\end{equation}
as well as
\begin{equation}
\label{ExandGx}
    E_h^x=E_h^{x,v}\text{ and }G_h^x=\bigcup_{v\in{V^0}}\big(E_h^{x,v}\cap F_h^{x,v}\big),
\end{equation}
and the (random) subsets of $\Z^d$
\begin{equation}
\label{EandG}
     E_h=\{x\in{\Z^d};\,E_h^x\text{ occurs}\}\text{ and }G_h=\{x\in{\Z^d};\,G_h^x\text{ occurs}\}.
\end{equation}
For all $K\subset\tilde{\Z}^d$, we denote by $\mathcal{A}_K$ the $\sigma$-algebra $\sigma\left(\tp_z,\ z\in{K}\right).$ We note that the sets $U^{x}$ are disjoint when $x$ varies, cf.\ \eqref{eq:Usets} and \eqref{eq:I_e}, and that the events $E_{h}^{x,v}$ are $\mathcal{A}_{\mathcal{K}^x}$-measurable. Theorem \ref{percphitilde} implies that $G_h$ contains an infinite connected component, and the goal is to go from this to the percolation of $E_h\cap\{x\in{\Z^d};\, \Phi_x\geqslant h\}.$ 
The following lemma makes the above observation, see the discussion at the beginning of this section, precise.
\begin{Lemme}
\label{comparingEandF}
There exists $h_1\in{(0,1]}$ such that for all $h\in(0,h_1]$ and $x\in{\Z^d},$
\begin{equation}
\label{comparingEandFeq}
\tilde{\P}^G\left(G_{h}^{x}\,\middle|\,\mathcal{A}_{\mathcal{K}^x}\right)\leqslant \tilde{\P}^G \left(E_h^x\cap\{\Phi_x\geqslant h\}\,\middle|\,\mathcal{A}_{\mathcal{K}^x}\right).
\end{equation}
\end{Lemme}

\begin{proof}
Let us fix some $x\in{\Z^d}.$  It is sufficient to prove that there exists $h_1\in{(0,1]}$ such that for all $h\in{(0,h_1]}$ and all $v \in V^0$,
\begin{equation}
\label{interm}
1_{E_h^{x,v}}\tilde{\P}^G\left(F_{h}^{x,v}\cap\{ \Phi_x\leqslant2h\}\,\middle|\,\mathcal{A}_{\mathcal{K}^x}\right)\leqslant\frac{1}{2d}1_{E_h^{x,v}}\tilde{\P}^G\left(h\leqslant\Phi_x\leqslant2h\,\middle|\,\mathcal{A}_{\mathcal{K}^x}\right).
\end{equation}
Indeed, if \eqref{interm} holds, then
\begin{align*}
\tilde{\P}^G\left(G_h^x\,\middle|\,\mathcal{A}_{\mathcal{K}^x}\right)&=\tilde{\P}^G\left(G_h^x\cap\{\Phi_x>2h\}\,\middle|\,\mathcal{A}_{\mathcal{K}^x}\right)+\tilde{\P}^G\left(G_h^x\cap\{\Phi_x\leqslant2h\}\,\middle|\,\mathcal{A}_{\mathcal{K}^x}\right)
\\[0.2em]&\stackrel{\eqref{ExandGx}}{\leqslant}\tilde{\P}^G\left(G_h^x\cap\{\Phi_x>2h\}\,\middle|\,\mathcal{A}_{\mathcal{K}^x}\right)+\sum_{v\in{V^0}}\tilde{\P}^G\left(E_h^{x,v}\cap F_h^{x,v}\cap\{\Phi_x\leqslant2h\}\,\middle|\,\mathcal{A}_{\mathcal{K}^x}\right)
\\[-0.2em]&\stackrel{\eqref{interm}}{\leqslant}\tilde{\P}^G\left(G_h^x\cap\{\Phi_x>2h\}\,\middle|\,\mathcal{A}_{\mathcal{K}^x}\right)+\frac{1}{2d}\sum_{v\in{V^0}}\tilde{\P}^G\left(E_{h}^{x,v}\cap\{h\leqslant\Phi_x\leqslant 2h\}\,\middle|\,\mathcal{A}_{\mathcal{K}^x}\right)
\\&\ \leqslant\tilde{\P}^G\left(E_h^x\cap\{\Phi_x\geqslant h\}\,\middle|\,\mathcal{A}_{\mathcal{K}^x}\right),
\end{align*}
noting that $G_h^x, E_h^{x,v}\subset E_h^x$ in the last inequality, and \eqref{comparingEandFeq} follows. We now show \eqref{interm}. Let us fix some $v\in{V^0}.$ We begin with the study of $\Phi_x,$ by decomposing it suitably. It follows from the Markov property, cf.\@ \eqref{phi^u}, that $\tilde{\Phi}_x^{U^x}=\Phi_x-\tilde{\beta}_x^{U^x}$ is a centered Gaussian variable with variance $g_{U^x}(x,x).$  The value of the variance $g_{U^x}(x,x)\equiv \sigma_0^2$ does not depend on $x\in{\Z^d}$ (it actually follows from Section 2 of \cite{MR3502602} that $\sigma_0^2=\frac{1}{4d}$). Moreover, on the event $E_h^{x,v},$ it is clear that $|\tilde{\beta}_x^{U^{x}}|\leqslant K(h).$ Thus, on the event $E_h^{x,v},$ since the harmonic average $\tilde{\beta}_x^{U^{x}}$ is $\mathcal{A}_{\mathcal{K}^x}$-measurable, we obtain, for all $h >0$,
\begin{align}
\label{firstineq}
    \tilde{\P}^G\left(-h\leqslant\Phi_x\leqslant2h\,|\,\mathcal{A}_{\mathcal{K}^x}\right)&=\tilde{\P}^G\left(-h\leqslant\tilde{\Phi}_x^{U^{x}}+\tilde{\beta}_x^{U^x}\leqslant2h\, \Big|\,\mathcal{A}_{\mathcal{K}^x}\right)\nonumber
    \\&= \frac{1}{\sqrt{2\pi \sigma_0^2}}\int_{-h}^{2h}{\exp\left(-\frac{\big(y-\tilde{\beta}_x^{U^x}\big)^2}{2\sigma_0^2}\right)\mathrm{d}y}\nonumber
    \\&= \frac{1}{\sqrt{2\pi \sigma_0^2}}\exp\left(-\frac{\big(\tilde{\beta}_x^{U^x}\big)^2}{2\sigma_0^2}\right)\int_{-h}^{2h}{\exp\left(-\frac{y^2}{2\sigma_0^2}\right)\exp\left(\frac{y\tilde{\beta}_x^{U^x}}{\sigma_0^2}\right)\mathrm{d}y}\nonumber
    \\&\leqslant \frac{1}{\sqrt{2\pi \sigma_0^2}}\exp\left(-\frac{\big(\tilde{\beta}_x^{U^x}\big)^2}{2\sigma_0^2}\right)\times3h\exp\left(\frac{2hK(h)}{\sigma_0^2}\right).
\end{align}
A similar calculation shows that on the event $E_h^{x,v},$ for $h >0$,
\begin{equation}
\label{secondineq}
     \tilde{\P}^G\left(h\leqslant\Phi_x\leqslant2h\,|\,\mathcal{A}_{\mathcal{K}^x}\right)\geqslant\frac{1}{\sqrt{2\pi \sigma_0^2}}\exp\left(-\frac{\big(\tilde{\beta}_x^{U^x}\big)^2}{2\sigma_0^2}\right)\times h\exp\left(-\frac{2h(K(h)+h)}{\sigma_0^2}\right).
\end{equation}
Define
\begin{equation*}
    C_6(d)=\sup_{h\in{(0,1]}}\left\{ 3\exp\left(\frac{2hK(h)}{\sigma_0^2}\right)\times\exp\left(\frac{2h(K(h)+h)}{\sigma_0^2}\right) \right\},
\end{equation*}
and note that $C_6<\infty$ since $hK(h)\rightarrow0$ as $h\searrow0,$ cf.\ \eqref{K}. Hence, by \eqref{firstineq} and \eqref{secondineq},
\begin{equation}
\label{inequalityforthegaussianpart}
     1_{E_h^{x,v}}\tilde{\P}^G\left(-h\leqslant\Phi_x\leqslant2h\,|\,\mathcal{A}_{\mathcal{K}^x}\right)\leqslant C_61_{E_h^{x,v}}\tilde{\P}^G\left(h\leqslant\Phi_x\leqslant2h\,|\,\mathcal{A}_{\mathcal{K}^x}\right), \text{ for }h\in{(0,1]}.
\end{equation}
Let us now turn to the events $F_h^{x,v}.$ It follows again from the Markov property for the Gaussian free field, see in particular the discussion below \eqref{bridgeinsideedge}, that, knowing $\mathcal{A}_{\mathcal{K}^x\cup\{x\}},$ the process $(\tilde{\Phi}_{z})_{z\in{U^{x,v}}}$ is a Brownian bridge of length $\frac{1}{4}$ between $\Phi_x$ and $\tilde{\Phi}_{x+\frac{1}{4}v}$ of a Brownian motion with variance $2$ at time $1.$ Using Lemma \ref{Brownianbridge}, one can then find $h_1\in{(0,1]}$ such that, for all $h\in{(0,h_1]},$ on the $\mathcal{A}_{\mathcal{K}^x\cup\{x\}}$ measurable event $E_h^{x,v}\cap\{-h\leqslant\Phi_x\leqslant2h\},$
\begin{align}
\label{inequalityforthebrigdepart}
\tilde{\P}^G\left(F_h^{x,v}\,\middle|\,\mathcal{A}_{\mathcal{K}^x\cup\{x\}}\right)&=\tilde{\P}^G\left(\min_{z\in{U^{x,v}}} \tilde{\Phi}_{z}\geqslant-h\,\middle|\,\mathcal{A}_{\mathcal{K}^x\cup\{x\}}\right)\nonumber
\\&=1-\exp\left(-4(h+\Phi_x)(h+\tilde{\Phi}_{x+\frac{1}{4}v})\right)\nonumber
\\&\leqslant1-\exp\left(-12h(K(h)+h)\right)\nonumber
\\&\leqslant\frac{1}{2dC_6}.
\end{align}
We now conclude using (\ref{inequalityforthegaussianpart}) and (\ref{inequalityforthebrigdepart}): for all $h\in(0,h_1],$
\begin{align*} 
1_{E_h^{x,v}}\tilde{\P}^G\left(F_{h}^{x,v}\cap\{\Phi_x\leqslant2h\}\,\middle|\,\mathcal{A}_{\mathcal{K}^x}\right)
&\stackrel{\eqref{FFFF}}{=}\E\left[1_{E_h^{x,v}\cap\{-h\leqslant\Phi_x\leqslant 2h\}}\tilde{\P}^G\left(F_h^{x,v}\,\middle|\,\mathcal{A}_{\mathcal{K}^x\cup\{x\}}\right)\,\middle|\,\mathcal{A}_{\mathcal{K}^x}\right]
\\&\stackrel{\eqref{inequalityforthebrigdepart}}{\leqslant}\frac{1}{2dC_6} 1_{E_h^{x,v}}\tilde{\P}^G\left(-h\leqslant\Phi_x\leqslant2h\,|\,\mathcal{A}_{\mathcal{K}^x}\right)
\\&\stackrel{\eqref{inequalityforthegaussianpart}}{\leqslant}\frac{1}{2d}1_{E_h^{x,v}}\tilde{\P}^G\left(h\leqslant\Phi_x\leqslant 2h\,\middle|\,\mathcal{A}_{\mathcal{K}^x}\right),
\end{align*}
which is \eqref{interm}.
\end{proof}
Lemma \ref{comparingEandF} roughly asserts that it is more likely to have $\{\Phi_x\geqslant h\}$ than to have $G_h^x$ (on $E_h^x$) for small $h >0$ and we know by Theorem \ref{percphitilde} that $G_h$ has an infinite connected component in a thick slab. Using the Markov property for the Gaussian free field we will show that this implies that $E^{\geqslant h},$ see \eqref{exc}, percolates in a sufficiently thick slab for any such value of $h$, thus obtaining Theorem \ref{mainresult}.

\begin{proof}[Proof of Theorem \ref{mainresult}]
Let us fix some $h\in{(0,h_1]},$ with $h_1$ as in Lemma \ref{comparingEandF}. The Markov property for the Gaussian free field, see \eqref{markov}, \eqref{eq:Usets} and \eqref{eq:I_e}, implies that the family $(\tilde{\Phi}_{\cdot}^{U^x})_{x\in{\Z^d}}$ is i.i.d.\  and independent of $\mathcal{A}_{\mathcal{K}},$ and that for all $x\in{\Z^d}$ and $z\in{U^x},$
\begin{equation*}
    \tilde{\Phi}_z=\tilde{\beta}^{U^x}_z+\tilde{\Phi}_z^{U^x}
\end{equation*} 
where $\tilde{\beta}^{U^x}$ is $\mathcal{A}_{\mathcal{K}^x}$-measurable for all $x\in{\Z^d}$. For each $x\in{\Z^d},$ there exists $f_x:C(\tilde{\Z}^d,\R)\times\R^{\mathcal{K}^x}\rightarrow\{0,1\}$ and $g_x:C(\tilde{\Z}^d,\R)\times\R^{\mathcal{K}^x}\rightarrow\{0,1\}$ such that
\begin{equation}
\label{fandg}
    1_{G_h^x}=f_x(\tilde{\Phi}^{U^x},\tilde{\Phi}_{|\mathcal{K}^x})\text{ and }1_{E_h^x\cap\{\Phi_x\geqslant h\}}=g_x(\tilde{\Phi}^{U^x},\tilde{\Phi}_{|\mathcal{K}^x}).
\end{equation}
For each measurable subset $A$ of $\tilde{\Z}^d,$ let us denote by $\tilde{\P}^{G,A}$ the law of $\tilde{\Phi}_{|A}$ under $\tilde{\P}^G.$ Lemma \nolinebreak\ref{comparingEandF} now gives that for all $x\in{\Z^d}$ and for $\tilde{\P}^{G,\mathcal{K}^x}$-a.s.\ all $\beta^x\in{\R^{\mathcal{K}^x}},$
\begin{equation}
\label{fx<gx}
    \tilde{\P}^G\big(f_x(\tilde{\Phi}^{U^x},\beta^x)=1\big)\leqslant\tilde{\P}^G\big(g_x(\tilde{\Phi}^{U^x},\beta^x)=1\big).
\end{equation}
For each $x\in{\Z^d}$ and $\beta^x\in{\R^{\mathcal{K}^x}}$ such that \eqref{fx<gx} holds, abbreviating the left and right-hand sides of \eqref{fx<gx} by $p_f(\beta^x)$ and $p_g(\beta^x)$, respectively, so that $p_f(\beta^x) \leqslant p_g(\beta^x)$, we can now define a probability $\nu_{\beta^x}$ on $\{0,1\}^2$ such that, with $\pi_1$ and $\pi_2$ respectively denoting the projections onto the first and second coordinate of $\{0,1\}^2,$ we have
\begin{equation}
\label{pi1and2}
 \pi_1\leqslant \pi_2, \ \   \nu_{\beta^x}(\pi_1=1)=p_f(\beta^x)\quad\text{ and }\quad  \nu_{\beta^x}(\pi_2=1)=p_g(\beta^x).
\end{equation}
The measure $\nu_{\beta^x}$ can for instance be constructed from a uniform random variable $Y$ on $[0,1]$ as the law of $(1_{\{ Y \leqslant p_f(\beta^x) \}}, 1_{\{ Y \leqslant p_g(\beta^x) \}})$ on $\{0,1\}^2$.
For each $\beta\in{\R^{\mathcal{K}}}$ and $x\in{\Z^d},$ let $\beta^x=(\beta_z)_{z\in{\mathcal{K}^x}}\in{\R^{\mathcal{K}^x}},$ and finally define, for $\tilde{\P}^{G,\mathcal{K}}$-a.s.\ all $\beta\in{\R^{\mathcal{K}}}$ the following probabilities on $(\{0,1\}^{\Z^d})^2$
\begin{equation}
\label{nubeta}
    \nu_{\beta}=\bigotimes_{x\in{\Z^d}}\nu_{\beta^x}\qquad\text{ and }\qquad\nu=\tilde{\E}^G\left[\nu_{\tilde{\Phi}_{|\mathcal{K}}}\right].
\end{equation}
Note that for all $B\subset\{0,1\}^2,$ $\beta^x\mapsto\nu_{\beta^x}(B)$ is measurable, and thus $\nu$ is well-defined. Let $\pi'_1$ and $\pi'_2$ be the projections on the first and second coordinate of $(\{0,1\}^{\Z^d})^2$. Then $\pi'_1\leqslant\pi'_2$ $\nu$-a.s by \eqref{pi1and2}, \eqref{nubeta}, and on account of \eqref{fandg} $\pi'_1$ has the same law under $\nu$ as $(1_{G_h^x})_{x\in{\Z^d}}$ under $\tilde{\P}^G$ and $\pi'_2$ has the same law under $\nu$ as $(1_{E_h^x\cap\{\Phi_x\geqslant h\}})_{x\in{\Z^d}}$ under $\tilde{\P}^G.$ Moreover, by Theorem \ref{percphitilde}, there exists $L_0=L_0(h)>0$ such that the set $\tilde{A}_h(\tilde{\Phi}),$ cf.\ \eqref{A''h}, contains $\tilde{\P}^G$-a.s.\ an unbounded connected component $M$ in the thick slab $\tilde{\Z}^2\times[0,2L_0)^{d-2}.$ By definition, see \eqref{EandG}, $M\cap\Z^d\subset G_h,$ and thus $\{x\in{\Z^d}: \pi'_1(x)=1\}$ contains $\nu$-a.s.\ an unbounded component in $\tilde{\Z}^2\times[0,2L_0)^{d-2}.$ Since $\pi'_1\leqslant \pi'_2,$ this implies that $E_h\cap \{x\in{\Z^d}:\,\tilde{\Phi}_x\geqslant h\}$ also contains an infinite connected component in $\tilde{\Z}^2\times[0,2L_0)^{d-2}$, as desired.
\end{proof}

\begin{Rk} \label{R:final}1) The result of \cite{MR3053773} is actually slightly better than Theorem \ref{mainresult} in high dimensions: if $d$ is large enough, there exist $h_2=h_2(d)>0$ and $L_0=L_0(d)\geqslant1$ such that the level set $\{x\in{\Z^d}; \, \Phi_x\geqslant h_2\}$ percolates in the slab $\Z^2\times[0,2L_0)\times\{0\}^{d-3}.$ However, in all dimensions $d\geqslant3,$ the set $\{x\in{\Z^d}; \Phi_x\geqslant h\}$ never percolates for $h\geqslant0$ in $\Z^2\times\{0\}^{d-2},$ as explained in Remark 3.6.1 of \cite{MR3053773}. 
\medskip

\noindent  2) It is possible to get a result similar to \eqref{1:30} for the positive level set of the Gaussian free field $\Phi$ on $\Z^d$ just constructed, thus obtaining the following strengthening of Theorem \ref{mainresult}. For all $h\leqslant h_1,$ let $K(h)$ be as in $\eqref{K}$ for $h_0=1,$ then the set $\{x\in{\Z^d};\,h\leqslant\Phi_x\leqslant K(h)\}$ contains a.s.\, an infinite connected component. Indeed, using an argument similar to that of Lemma \ref{comparingEandF}, one can prove that, conditionally on $\mathcal{A}_{\mathcal{K}^x},$ the probability of $G_h^x\cap\{\Phi_x\leqslant K(h)\}$ is smaller than the probability of $E_h^x\cap\{h\leqslant\Phi_x\leqslant K(h)\}$ and the result follows.
\medskip

\noindent 3) Theorem 2.2 in \cite{MR3024098} can also easily be extended to the Gaussian free field: for each $h\leqslant h_1,$ the set $\{x\in{\Z^d};\,\Phi_x\geqslant h\}$ contains an almost surely transient component. Indeed, looking at the proof of Theorem 2.2 in \cite{MR3024098}, see also Theorem 1 in \cite{MR2819660}, we can use \eqref{H0N} instead of (5.1) in \cite{MR3024098} to obtain that the set $\tilde{A}'_{u,p}$ defined in \eqref{A'u} contains an unbounded connected and transient component for $p\in{[p(u),1]}.$ Using the same coupling as in Lemma \ref{couplingwithbernoulli}, we get that this is also true for the set $\tilde{A}_h(\tilde{\phi})$ defined in \eqref{A''h}, and the same proof as the proof of Theorem \ref{mainresult} tells us that $\{x\in{\Z^d}; \, \Phi_x\geqslant h\}$ also contains an infinite connected and transient component for $h\leqslant h_1.$
\medskip

\noindent 4) Another parameter $\overline{h}\leqslant h_*$ has been introduced in \cite{MR3390739}, and a similar one has been used in \cite{MR3417515}. This parameter describes a strong percolative regime for $E^{\geqslant h}$, when $h<\overline{h},$ i.e., all connected components of $E^{\geqslant h}$ in $[-R,R]^d$ with diameter at least $\frac{R}{10}$ are connected in $[-2R,2R]^d$ with large enough probability when $R$ goes to $\infty.$ It has been proved that $\overline{h}>-\infty$ and it is believed that actually $\overline{h}=h_*,$ but it is still unknown whether $\overline{h}\geqslant0$ or not. 
Our methods may perhaps help in that regard.
\end{Rk}

\noindent\textbf{Acknowledgments.}
We thank the anonymous referees for their careful reading of the manuscript and for their numerous comments. AD gratefully acknowledges support of the 
UoC Forum \lq Classical and quantum dynamics of interacting particle systems\rq, and AP gratefully acknowledges support of the \lq IPaK\rq-program at the University of Cologne.

\appendix
\setcounter{secnumdepth}{0}
\section{Appendix: Proof of Lemma \ref{connectivity}}
\renewcommand*{\theThe}{A.\arabic{The}}
\setcounter{The}{0}
\setcounter{equation}{0}
\renewcommand{\theequation}{A.\arabic{equation}}

The proof of Lemma \ref{connectivity} is very close to the proof of Proposition 1 in \cite{MR2819660}, but we need to remove the dependence on $u$ of the constants, and make the dependence on $u$ of the error term explicit instead. 
We will henceforth refer to \cite{MR2819660} whenever possible, and in particular, 
Lemmas 3 to 6 and 11 in \cite{MR2819660} do not involve $u$ at all, so we will use them without proof. 
Recall $\omega^u$, the interlacement process at level $u$ on $\Z^d,$ and $\tilde{\omega}^u$, the interlacement process at level $u$ on the cable system, obtained from $\omega^u$ by adding independent Brownian excursions on the edges, as in the construction of the diffusion $\tilde{X}$ from a simple random walk on $\Z^d$ in the beginning of section \ref{cablesystem}. We denote by $\hat{\I}^u$ the set of edges traversed by at least one of the trajectories in $\text{supp}(\omega^u).$ Now observe that the event that every $x$ and $y$ in $\tilde{\I}^u\cap[0,R)^d$ be connected in $\tilde{\I}^u\cap[-\epsilon R,(1+\epsilon)R)^d,$ which is the event of interest in \eqref{eq:connectivity}, is more likely than every $x$ and $y$ in $\hat{\I}^u\cap[-1,R+1)^d$ being connected in $\hat{\I}^u\cap[-\epsilon R,(1+\epsilon)R)^d.$ Thus, we only need to show the respective statement of Lemma \ref{connectivity} for $\hat{\I}^u$ instead of $\tilde{\I}^u$, cf.\ Lemma \ref{13} below. 

The idea of the proof is to show that there exists $C\geqslant1$ such that for every integer $R\geqslant1$ and every $x,y\in{\I^u\cap[-R,R)^d},$ the vertices $x$ and $y$ are connected through edges in $\hat{\I}^u\cap[-CR,CR)^d$ with high enough probability. It is quite hard to directly link $x$ and $y,$ especially if $R$ is large. Therefore, let us define, under some probability $P,$ $\omega_{i,3}^{u/3}$ for $i\in{\{1,2,3\}}$, three independent Poisson point process with the same law as $\omega^{u/3}$ under $\P^I,$ such that $\omega^{u}=\sum_{i=1}^3\omega_{i,3}^{u/3}.$ Let us call $\I_{i,3}^{u/3}$ the set of vertices visited by at least one of the trajectories from supp$(\omega_{i,3}^{u/3})$, denote by $\hat{\I}_{i,3}^{u/3}$ the set of edges traversed by at least one of the trajectories from supp$(\omega_{i,3}^{u/3}),$ and let $C_i^{u/3}(x,R)$ be the set of vertices connected to $x$ by edges in $\hat{\I}^{u/3}_{i,3}\cap[-R,R)^d$ for $i\in{\{1,2,3\}}$ and $x\in{\Z^d}.$ We are going to prove that, if $x\in{\I_{1,3}^{u/3}}$ and $y\in{\I_{2,3}^{u/3}},$ then $C_1^{u/3}(x,R)$ and $C_2^{u/3}(y,R)$ are big enough, and that one can connect these two sets by edges in $\hat{\I}_{3,3}^{u/3}\cap[-CR,CR)^d$ with high probability. In particular, this will imply that $x$ and $y$ are connected through edges in $\hat{\I}^u\cap[-CR,CR)^d$ with high probability.
\medskip

We first recall a property of the Poisson distribution (see for example (2.11) in \cite{MR2819660}): let $N$ be a random variable which has Poisson distribution with parameter $\lambda,$ then there exist constants $c<1$ and $C>1$ independent of $\lambda$ such that
\begin{equation*}
    \P\left(c\lambda\leqslant N\leqslant C\lambda\right)\geqslant1-C\exp\left(-c\lambda\right).
\end{equation*}
For $A\subset\Z^d$ finite and $\omega$ an interlacement process, we call for all $u>0$ $N_A^u$ the the number of trajectories in supp$(\omega^u)$ which enter $A$, and write $Z_1,\dots,Z_{N_A^u}$ for the corresponding trajectories, parametrized such that $Z_i(0)\in{A}$ and $Z_i(-n)\notin{A}$ for all $n>0.$ Note that $Z_1,\dots,Z_{N_A^u}$ depend on $\omega,$ $u$ and $A$ even if this is only implicit in the notation. Then $N_A^u$ is a Poisson variable with parameter $u\mathrm{cap}(A)$ and 
\begin{equation}
\label{poisson}
    \P^I\left(c u\mathrm{cap}(A)\leqslant N_A^u\leqslant Cu\mathrm{cap}(A)\right)\geqslant1-C\exp\left(-cu\mathrm{cap}(A)\right).
\end{equation}
Here, $\mathrm{cap}(A)$ is the capacity of the set $A$, i.e., the total mass of the equilibrium measure of $A$. The following standard bounds will soon prove to be useful: For any $A\subset[-R,R)^d\cap\Z^d$ and $R \geqslant 1$,
\begin{equation}
\label{capacity}
\mathrm{cap}(A)\leqslant\mathrm{cap}([-R,R)^d)\leqslant CR^{d-2}\text{ and }\mathrm{cap}([-R,R)^d)\geqslant cR^{d-2}.
\end{equation}
The next lemma gives a bound on the probability to connect the two sets $C_1^{u/3}(x,R)$ and $C_2^{u/3}(y,R)$ in $\hat{\I}_3^{u/3}\cap[-CR,CR)^d$ in terms of capacity.

\begin{Lemme}
\label{12}
There exist constants $c=c(d)>0$ and $C=C(d)<\infty$ such that for all $R>0$ and $u>0,$ for all subsets $U$ and $V$ of $[-R,R)^d,$ 
\begin{equation*}
    \P^I\big(U\stackrel{\hat{\I}^u\cap [-CR,CR)^d}\longleftrightarrow V\big)\geqslant1-C\exp\left(-cR^{2-d}u\mathrm{cap}(U)\mathrm{cap}(V)\right).
\end{equation*}
\end{Lemme}

\begin{proof}
If there is a trajectory among $(Z_1,\dots,Z_{N_U^u}),$ which hits $V$ after $0$ and before leaving $[-CR,CR)^d$, then $U$ is connected to $V$ through edges of $\hat{\I}^u\cap [-CR,CR)^d.$ We can use Lemma 11 in \cite{MR2819660} to lower bound the probability of a trajectory to behave accordingly by $cR^{2-d}\mathrm{cap}(V)$, and thus we infer
\begin{align*}
    \P^I\left(U\stackrel{\hat{\I}^u\cap[-CR,CR)^d}\longleftrightarrow V\right)&\geqslant 1-\P^I(N_U^u<cu\mathrm{cap}(U))-\left(1-cR^{2-d}\mathrm{cap}(V)\right)^{cu\mathrm{cap}(U)}
    \\&\stackrel{\mathclap{(\ref{poisson}),\,(\ref{capacity})}}\geqslant\phantom{.2)}1-C\exp\left(-cR^{2-d}u\mathrm{cap}(U)\mathrm{cap}(V)\right).
\end{align*}
$\quad$
\end{proof}
We are now going to prove that $\mathrm{cap}\big(C_1^{u/3}(x,R)\big)$ and $\mathrm{cap}\big(C_2^{u/3}(y,R)\big)$ are large enough with high probability, and in particular that they grow faster in $R$ than $R^{\frac{d-2}{2}}.$ From now on we fix some $u_0>0,$ and for all $u>0,$ $A \subset \Z^d$ finite and $T$ a positive integer, we define the set $\Psi(u,A,T)$ by
\begin{equation}
\label{Psi}
    \Psi(u,A,T)=A\cup\bigcup_{i=1}^{N_A^u}{\{Z_i(n),\ 0\leqslant n\leqslant T\}}.
\end{equation}

\begin{Lemme}
\label{firstboundtry}
For all $\eps\in{(0,1)},$ $k\geqslant1$ and $\delta\geqslant\eps,$ there exist constants $c>0$ and $C<\infty$ such that for every $u\in{(0,u_0]}$,  $A\subset\Z^d$ finite and $T$ a positive integer, 
\begin{equation}
\label{boundoncappsi}
\begin{split}
&\P^I\left(\mathrm{cap}\left(\Psi\left(u,A,T\right)\right)\geqslant c\, \min\left(u\mathrm{cap}(A)T^{\frac{1-\epsilon}{2}},T^{\frac{(d-2)(1-\epsilon)}{2}}\right)\right)\\
&\qquad \qquad\qquad \qquad \qquad \qquad \geqslant 1-C\exp\left(-c\min\left(T^{\epsilon/2},u\mathrm{cap}(A)\right)\right),
\end{split}
\end{equation}
and, if $A\subset B=[-kT^\delta,kT^\delta)^d$,
\begin{equation}
\label{boundonpsi}
\P^I\left(\Psi(u,A,T)\subset B+[-T^{\frac{1+\eps}{2}},T^{\frac{1+\eps}{2}})^d\right)\geqslant 1-C\exp\left(-cT^\eps u\right).
\end{equation}
\end{Lemme}

\begin{proof}
(\ref{boundoncappsi}) is a simple consequence of Lemma 6 in \cite{MR2819660} and (\ref{poisson}). In order to prove (\ref{boundonpsi}), let us first define $h(T,\eps)$, the probability that the simple random walk on $\Z^d$ beginning in $0$ leaves $[-T^{\frac{1+\eps}{2}},T^{\frac{1+\eps}{2}})^d$ before time $T$. Hoeffding's inequality yields that $h(T,\eps)\leqslant C\exp(-cT^\eps).$ Now, taking $B=[-kT^\delta,kT^\delta)^d$ and assuming $A\subset B$, observing that, in order for $\Psi(u,A,T)$ not to be contained in $ B+[-T^{\frac{1+\eps}{2}},T^{\frac{1+\eps}{2}})^d$, at least one of the walks $Z_i$ in \eqref{Psi} must reach distance $T^{\frac{1+\eps}{2}}$ before time $T$, and noting that $N_A^u \leq N_B^u$, we get 
\begin{align*}
    \P^I\left(\Psi(u,A,T)\subset B+[-T^{\frac{1+\eps}{2}},T^{\frac{1+\eps}{2}})^d\right)&\geqslant1-\P^I\left(N_B^u\geqslant Cu\mathrm{cap}(B)\right)-Cu\mathrm{cap}(B)h(T,\eps)
    \\&\stackrel{\mathclap{(\ref{poisson}),\,(\ref{capacity})}}\geqslant1-C\exp\left(-cuT^{(d-2)\delta}\right)-CuT^{(d-2)\delta}\exp\left(-cT^\eps\right)
    \\[0.2em]&\geqslant1-C\exp\left(-cT^\eps u\right),
\end{align*}
which concludes the proof.
\end{proof}

We now iterate this process to find the desired bound on $\mathrm{cap}\left(C^{u}(x,R)\right).$ Consider, under some probability $\Q,$ a sequence of independent random interlacement processes $(\omega_k)_{k\geqslant1}$ which define an independent sequence $\left(\Psi_k\right)_{k\geqslant2}$ such that for all $k\geqslant2,$ $\Psi_k$ has the same law as $\Psi$ (see (\ref{Psi}) for notation), and let $\mathcal{I}_1^v$ be the random interlacement set associated with $\omega_1^v.$ For each $x\in{\Z^d},$ let $Z^x$ be the trajectory with the smallest label $v$ contained in $\omega_1$ such that $Z^x(0)=x,$ which exists since $x\in{\mathcal{I}_1^v}$ for $v\in{(0,\infty)}$ large enough. For all $x\in{\Z^d},$ $u>0$ and $T$ positive integer, we recursively define a sequence of subsets  $(U^{(k)}_u(x,T))_{k\geqslant 1}$ of $\Z^d$ by 
\begin{equation*}
    U_u^{(1)}(x,T)=\{Z^x(n),\ 0\leqslant n\leqslant T\},
\end{equation*}
and, for all $k\geqslant2,$
\begin{equation*}
    U_u^{(k)}(x,T)=\Psi_k\left(u,U_u^{(k-1)}(x,T),T\right).
\end{equation*}
Note that $U_u^{(k)}(x,T)$ depends on $u$ only for $k\geqslant2.$ In the next lemma, we iterate the results of Lemma \ref{firstboundtry} to find lower bounds on the capacity of $U_u^{(d-2)}(x,T)$ and upper bounds on the diameter of $U_u^{(d-2)}(x,T).$

\begin{Lemme}
\label{secondboundtry}
For all $\eps\in{(0,\frac{1}{3}]},$ there exist constants $c>0$ and $C<\infty$ such that for every $u\in{(0,u_0]},$ $x\in{\Z^d},$ positive integer T and $k\in{\{1,\dots,d-2\}},$ 
\begin{equation}
\label{4.12}
\Q\left(\mathrm{cap}\left(U_u^{(k)}(x,T)\right)
\geqslant u^{k-1}\big(cT^{\frac{1-\eps}{2}}\big)^k\right)
\geqslant 1-C\exp\left(-cT^{\epsilon/2}u\right)
\end{equation}
and
\begin{equation}
\label{4.11}
\Q\left(U_u^{(k)}(x,T)\subseteq x+[-kT^{\frac{1+\eps}{2}},kT^{\frac{1+\eps}{2}})^d\right)\geqslant 1-C\exp\left(-cT^{\eps}u\right).
\end{equation}
\end{Lemme}
\begin{proof}
Let us introduce the shorthand $U_u^{(k)}=U_u^{(k)}(x,T),$ and note that, albeit only implicitly, $U_u^{(k)}$ depends on $x$ and $T$ . We first prove (\ref{4.12}) by induction on $k\in{\{1,\dots,d-2\}}.$ For $k=1,$ \eqref{4.12} follows directly from Lemma 6 in \cite{MR2819660}. Let us assume that \eqref{4.12} holds true at level $k-1$ for some $k\in{\{2,\dots,d-2\}},$ and that $c$ is small enough so that $u_0c\leqslant1.$ The event in (\ref{4.12}) is implied by the event
\begin{equation*}
    \left\{\mathrm{cap}\left(U_u^{(k-1)}\right)
\geqslant u^{k-2}\big(cT^{\frac{1-\eps}{2}}\big)^{k-1}\right\}\cap\left\{\mathrm{cap}\left(U_u^{(k)}\right)
\geqslant c\min\left(u\mathrm{cap} \left(U_u^{(k-1)}\right)T^{\frac{1-\epsilon}{2}},T^{\frac{(d-2)(1-\epsilon)}{2}}\right)\right\}
\end{equation*}
since $k\leqslant d-2.$ We only need to prove (\ref{4.12}) if $cT^{\eps/2}u\geqslant1,$ and then $cT^{\frac{1-\eps}{2}}u\geqslant T^{\eps/2}$ and (\ref{boundoncappsi}) gives that, on the event $\mathrm{cap}\left(U_u^{(k-1)}\right)\geqslant u^{k-2}\left(cT^{\frac{1-\eps}{2}}\right)^{k-1},$
\begin{IEEEeqnarray*}{rcl}
    \IEEEeqnarraymulticol{3}{l}{\Q\left(\mathrm{cap}\left(U_u^{(k)}\right)
\geqslant c\min\left(u\mathrm{cap} \left(U_u^{(k-1)}\right)T^{\frac{1-\epsilon}{2}},T^{\frac{(d-2)(1-\epsilon)}{2}}\right)\,\middle|\,U_u^{(k-1)}\right)}
\\\qquad&\geqslant&1-C\exp\left(-c\min\left(T^{\epsilon/2},\left(cT^{\frac{1-\eps}{2}}u\right)^{k-1}\right)\right) \geqslant1-C\exp\left(-cT^{\epsilon/2}u\right),
\end{IEEEeqnarray*}
and \eqref{4.12} follows by the induction hypothesis. The proof of (\ref{4.11}) is similar: one needs to use the fact that $h(T,\eps)\leqslant C\exp(-cT^\eps)$ for $k=1$ (see the proof of Lemma \ref{firstboundtry} for the definition of $h(T,\eps)$), and then proceed by induction with (\ref{boundonpsi}) for $k\geqslant2.$
\end{proof}

\begin{Cor}
\label{10}
For all $\eps\in{(0,\frac{1}{2}]},$ there exist constants $c>0$ and $C<\infty$ such that for every $u\in{(0,u_0]},$ $x\in{\Z^d}$ and $R>0$, 
\begin{equation*}
\P^I\left(x\in{\I^u},\ \mathrm{cap}(C^u(x,R))<cR^{(1-\eps)(d-2)}u^{d-3}\right)\leqslant C\exp\left(-cR^{\epsilon/2}u\right),
\end{equation*}
where $C^u(x,R)$ is the set of vertices connected to $x$ by edges in $\hat{\I}^u\cap[-R,R)^d.$
\end{Cor}

\begin{proof}
For all $v\in{(0,\frac{u_0}{d-2}]},$ let $\omega^{(d-2)v}:=\sum_{i=1}^{d-2}\omega_i^{v},$ where $(\omega_i^{v})_{i\geqslant1}$ are the independent random interlacement processes at level $v$ used in the definition of $\big(U_{v}^{(k)}(.\,,.)\big)_{k\geqslant1},$ see above Lemma \ref{secondboundtry}, and let $\I^{(d-2)v}$ be the random interlacement set associated with $\omega^{(d-2)v}.$ By definition, $\omega^{(d-2)v}$ has the same law under $\Q$ as a random interlacement process at level $(d-2)v,$ and if $x\in{\I^{v}_1}$ then $U_{v}^{(k)}(x,T)$ is a connected subset of $\I^{(d-2)v}$ for all $k\in{\{1,\dots,d-2\}},$  $x\in{\Z^d}$ and positive integer $T.$ In particular, if $x\in{\I^{v}_1}$ and if the event in \eqref{4.11} occurs with $\varepsilon$ in that formula taking the value of some $\delta\in{(0,\frac{1}{3})},$ then
\begin{equation*}
U_{v}^{(d-2)}(x,T)\subset C^{(d-2)v}\big(x,(d-2)T^{\frac{1+\delta}{2}}\big).
\end{equation*} 
Using Lemma \ref{secondboundtry} with $k=d-2$ we obtain, for all $\delta\in{(0,\frac{1}{3})},$ $v\in{(0,u_0]},$ $x\in{\Z^d}$ and positive integer $T$,
\begin{equation*}
\Q\left(x\in{\I^{v}_1},\ \mathrm{cap}\left(C^{(d-2)v}(x,(d-2)T^{\frac{1+\delta}{2}})\right)< cT^{\frac{(1-\delta)(d-2)}{2}}v^{d-3}\right)
\leqslant C\exp\left(-cT^{\delta/2}v\right).
\end{equation*}
The result follows by taking $v=\frac{u}{d-2},$ $T=\left\lfloor\left(\frac{R}{d-2}\right)^{2-\eps}\right\rfloor$ and $\delta=\frac{\eps}{2-\eps}.$
\end{proof}

We now have all the tools required to connect $x,y\in{\I^u\cap[-R,R)^d}$ through edges in $\hat{\I}^u\cap[-CR,CR)^d$, as mentioned at the beginning of the Appendix. 

\begin{Lemme}
\label{13}
There exist constants $c>0$ and $C<\infty$ such that for every $u\in{(0,u_0]},$ $R>0$ and $x,y\in{[-R,R)^d},$ 
\begin{equation*}
    \P^I\left(x,y\in\I^u,\left\{x\stackrel{\hat{\I}^u\cap[-CR,CR)^d}\longleftrightarrow y\right\}^c\right)
    \leqslant C\exp\left(-cR^{1/6}u\right).
\end{equation*}
\end{Lemme}

\begin{proof}
Using the notation introduced at the beginning of the Appendix, we have
\begin{equation*}
\P^I\left(x,y\in\I^u,\left\{x\stackrel{\hat{\I}^u\cap[-CR,CR)^d}\longleftrightarrow y\right\}^c\right)
    \leqslant\sum_{i,j=1}^3{P\left(x\in\I_{i,3}^{u/3},y\in\I_{j,3}^{u/3},\left\{x\stackrel{\hat{\I}^u\cap[-CR,CR)^d}\longleftrightarrow y\right\}^c\right)}.
\end{equation*}
Let us now fix $i,j\in{\{1,2,3\}}$ and let $k\in{\{1,2,3\}}$ be different from $i$ and $j.$ We define the events
\begin{equation*}
    E_1=\left\{\mathrm{cap}\left(C_i^{u/3}(x,R)\right)\geqslant CR^{\frac{2(d-2)}{3}}u^{d-3}\right\}, \  E_2=\left\{\mathrm{cap}\left(C_j^{u/3}(y,R)\right)\geqslant CR^{\frac{2(d-2)}{3}}u^{d-3}\right\},
\end{equation*}
and note that $E_1\subset\{x\in{\I_{i,3}^{u/3}}\}$ and $E_2\subset\{y\in{\I_{j,3}^{u/3}}\}.$ Thus,
\begin{align}
\label{proofofLemma13}
 &P\left(x\in\I_{i,3}^{u/3},y\in\I_{j,3}^{u/3},\left\{x\stackrel{\hat{\I}^u\cap[-CR,CR)^d}\longleftrightarrow y\right\}^c\right)\nonumber \\
 \begin{split}
    &\leqslant P\left((E_1\cap E_2)\setminus\Big\{C_i^u(x,R)\stackrel{\hat{\I}_{k,3}^{u/3}\cap [-CR,CR)^d}\longleftrightarrow C_j^u(y,R)\Big\}\right)\\
&\quad +P\left(\{x\in{\I_{i,3}^{u/3}}\}\setminus E_1\right)+P\left(\{y\in{\I_{j,3}^{u/3}}\}\setminus E_2\right).
\end{split}
\end{align}
For $uR^{1/6}\geqslant1,$ we can now use Lemma \ref{12} to bound the first summand of (\ref{proofofLemma13}) as
\begin{align*}
    P\bigg((E_1\cap E_2)\setminus\Big\{C_i^u(x,R)\stackrel{\hat{\I}_{k,3}^{u/3}\cap [-CR,CR)^d}\longleftrightarrow C_j^u(y,R)\Big\}\bigg)&\leqslant C\exp\left(-cR^{2-d}u\times R^{\frac{4(d-2)}{3}}u^{2(d-3)}\right)
    \\&\leqslant C\exp\left(-cR^{1/6}u\right).
\end{align*}
The second summand of (\ref{proofofLemma13}) can also be bounded using Corollary \ref{10} with $\eps=\frac{1}{3},$ and the result follows. 
\end{proof}

We now come to the 
\begin{proof}[Proof of Lemma \ref{connectivity}]
Lemma \ref{connectivity} is a simple consequence of Lemma \ref{13}. Indeed, let us define $R'=\lfloor \eps R/2C\rfloor$ with $C$ as in Lemma \ref{13}, and we can assume without loss of generality that $\eps R\geqslant2C.$ We define for each $x\in{\Z^d}$ the events
\vspace{-0.2cm}
\begin{equation*}
A_x^{(1)}=\left\{\I^u\cap\left(x+[-R',R')^d\right)\neq\emptyset\right\}\text{ and }A_x^{(2)}=\bigcap_{x,y\in{\I^u\cap (x+[-2R',2R')^d)}}\left\{x\stackrel{\hat{\I}^u\cap x+[-\eps R,\eps R)^d}\longleftrightarrow y\right\}.
\end{equation*}
Note that these events depend on our choice of $u$ and $R$ even if it does not appear in the notation. It follows from the definition of random interlacements, (\ref{capacity}) and Lemma \ref{13} that
\vspace{-0.2cm}
\begin{equation*}
    \P^I\left(A_x^{(1)}\right)\geqslant1-\exp\left(-cR^{d-2}u\right)\text{ and }\P^I\left(A_x^{(2)}\right)\geqslant1-CR^{2d}\exp\left(-cR^{1/6}u\right).
\end{equation*}
\vspace{-0.2cm}
In particular, we get that
\vspace{-0.1cm}
\begin{equation*}
    \P^I\bigg(\bigcap_{x\in{[0,R)^d\cap\Z^d}}A_x^{(1)}\cap A_x^{(2)}\bigg)\geqslant1-CR^{3d}\exp\left(-cR^{1/6}u\right)\geqslant1-C\exp\left(-cR^{1/7}u\right).
\end{equation*}
Let us call $A$ the event on the left-hand side in the previous line, and suppose that $A$ occurs. Then, for all $x,y\in{\I^u\cap[0,R)^d},$ one can find a path of nearest neighbors between $x$ and $y$ in $[0,R)^d.$ Moreover, if $x$ and $x'$ are two neighbors in $[0,R)^d,$ then $(x+[-R',R')^d)\cup(x'+[-R',R')^d)\subset x+[-2R',2R')^d,$ so every vertex in 
\begin{equation}
\label{set}
\I^u\cap (x+[-R',R')^d)\text{ is connected to every vertex in }\I^u\cap (x'+[-R',R')^d)
\end{equation}
by a path of edges in $\hat{\I}^u\cap (x+[-\eps R,\eps R)^d)\subset\hat{\I}^u\cap[-\eps R,(1+\eps)R)^d,$ and the sets in \eqref{set} are not empty. This tells us that if $A$ occurs, then every $x,y\in{\I^u\cap[0,R)^d}$ can be connected by edges in $\hat{\I}^u\cap[-\eps R,(1+\eps)R)^d,$ and thus $A$ implies the event on the left-hand side of \eqref{connectivity}.
\end{proof}

\bibliography{bibliographie}
\bibliographystyle{plain}
\bibliographystyle{short}

\end{document}